\documentclass[a4paper,12pt]{amsart}
\usepackage{amssymb,amsmath,amsthm}
\usepackage{ifthen}
\usepackage{graphicx}
\nonstopmode \numberwithin{equation}{section}
\newtheorem{thm}[equation]{Theorem}
\newtheorem{cor}[equation]{Corollary}
\newtheorem{lem}[equation]{Lemma}
\newtheorem{prop}[equation]{Proposition}

\newtheorem{conj}{Conjecture}

\theoremstyle{definition}

\newtheorem{prob}[equation]{Problem}

\newenvironment{rem}{%
\bigskip
\noindent \textsl{{\sl Remark. }}}{\bigskip}
\newenvironment{rems}{%
\bigskip
\noindent \textsl{{\sl Remarks. }}}{\bigskip}

\newcounter {own}
\def\theown {\thesection       .\arabic{own}}

\newenvironment{nonsec}{\bf
\setcounter{own}{\value{equation}}\addtocounter{equation}{1}
\refstepcounter{own}
\bigskip

\noindent \theown.$\ \,$}{$\,$.\ \ \ }

\newenvironment{pf}[1][]{%
 \vskip 3mm
 \noindent
 \ifthenelse{\equal{#1}{}}%
  {{\slshape Proof. }}%
  {{\slshape #1.} }%
 }%
{\qed\bigskip}

\newcounter{alphabet}
\newcounter{tmp}

\newcommand{\ID}{{\mathbb D}}

\newcommand{\IC}{{\mathbb C}}

\def\be{\begin{equation}}
\def\ee{\end{equation}}

\newcommand{\bee}{\begin{enumerate}}
\newcommand{\eee}{\end{enumerate}}

\newcommand{\blem}{\begin{lem}}
\newcommand{\elem}{\end{lem}}
\newcommand{\bthm}{\begin{thm}}
\newcommand{\ethm}{\end{thm}}
\newcommand{\bcor}{\begin{cor}}
\newcommand{\ecor}{\end{cor}}
\newcommand{\beg}{\begin{examp}}
\newcommand{\eeg}{\end{examp}}
\newcommand{\begs}{\begin{examples}}
\newcommand{\eegs}{\end{examples}}
\newcommand{\bdefe}{\begin{defin}}
\newcommand{\edefe}{\end{defin}}
\newcommand{\bprob}{\begin{prob}}
\newcommand{\eprob}{\end{prob}}
\newcommand{\bei}{\begin{itemize}}
\newcommand{\eei}{\end{itemize}}

\newcommand{\bcon}{\begin{conj}}
\newcommand{\econ}{\end{conj}}
\newcommand{\bcons}{\begin{conjs}}
\newcommand{\econs}{\end{conjs}}
\newcommand{\bprop}{\begin{prop}}
\newcommand{\eprop}{\end{prop}}
\newcommand{\br}{\begin{rem}}
\newcommand{\er}{\end{rem}}
\newcommand{\brs}{\begin{rems}}
\newcommand{\ers}{\end{rems}}
\newcommand{\bo}{\begin{obser}}
\newcommand{\eo}{\end{obser}}
\newcommand{\bos}{\begin{obsers}}
\newcommand{\eos}{\end{obsers}}
\newcommand{\bpf}{\begin{pf}}
\newcommand{\epf}{\end{pf}}
\newcommand{\ba}{\begin{array}}
\newcommand{\ea}{\end{array}}
\newcommand{\beq}{\begin{eqnarray}}
\newcommand{\beqq}{\begin{eqnarray*}}
\newcommand{\eeq}{\end{eqnarray}}
\newcommand{\eeqq}{\end{eqnarray*}}

\newcommand{\ds}{\displaystyle}

\begin{document}
\bibliographystyle{amsplain}
\title{Region of variability for certain classes of univalent functions
satisfying differential inequalities}
\author{S. Ponnusamy}
\address{S. Ponnusamy, Department of Mathematics,
Indian Institute of Technology Madras, Chennai-600 036, India.}
\email{samy@iitm.ac.in}
\author{A. Vasudevarao}
\address{A. Vasudevarao, Department of Mathematics,
Indian Institute of Technology Madras, Chennai-600 036, India.}
\email{alluvasu@iitm.ac.in}
\author{M. Vuorinen}
\address{M. Vuorinen, Department of Mathematics,
FIN-20014 University of Turku, Finland.}
\email{vuorinen@utu.fi}

\date{\today,
; File: pvdev$_{-}$7final.tex
}

\begin{abstract}
For complex numbers $\alpha$, $\beta$ and $M\in\mathbb{R}$ with $0<M\leq |\alpha|$
and $|\beta|\leq 1$,
let $\mathcal{B}(\alpha, \beta,M)$ be the class of analytic and univalent functions $f$ in
the unit disk ${\mathbb D}$ with $f(0)=0$, $f'(0)=\alpha$ and $f''(0)=M\beta$ satisfying
$|zf''(z)|\leq M$, $z\in\mathbb{D}$.
Let $\mathcal{P}(\alpha, M)$ be the another class of analytic and univalent functions in
$\mathbb{D}$ with $f(0)=0$, $f'(0)=\alpha$ satisfying ${\rm Re\,}(zf''(z))>-M$,
$z\in\mathbb{D}$,
where $\alpha\in\mathbb{C}\setminus\{0\}$, $0<M\leq 1/\log 4$.
For any fixed $z_0 \in {\mathbb D}$ and $\lambda\in\overline{\mathbb{D}}$
we shall determine the region of variability $V_j$ $(j=1,2)$ for $f'(z_0)$ when $f$
ranges over the class $\mathcal{S}_j$ $(j=1,2)$, where
\beqq
\mathcal{S}_1  & = &
\left\{\frac{}{} f\in\mathcal{B}(\alpha, \beta,M)
:\,f'''(0)=M(1-|\beta|^2)\lambda 
\right\}
\eeqq
and
$$\mathcal{S}_2 = \left \{f\in\mathcal{P}(\alpha, M)
:\,f''(0)=2M\lambda 
\right \}.
$$
In the final section we graphically illustrate the region of variability
for several sets of parameters.

\end{abstract}

\thanks{The first author was supported by NBHM (DAE, sanction No. 48/2/2006/R\&D-II),
while the second author was supported by
NBHM (DAE, sanction No. 48/2/2006/R\&D-II) and CIMO
(Grant no.15.5.2007/TM-07-5076/CIMO Fellowship),
Academy of Finland, and  Research project ``Quasiconformal Mappings'' of Matti Vuorinen}

\maketitle \pagestyle{myheadings} \markboth
{S. Ponnusamy, A. Vasudevarao and M. Vuorinen}{Regions of variability}

\section{Introduction and Preliminaries}
We denote the class of all analytic functions in the unit disk
$\ID=\{z\in\IC:\,|z|<1\}$ by ${\mathcal H}(\ID)$, and  think of
$\mathcal{H}(\ID)$ as a topological vector space endowed with the
topology of uniform convergence over compact subsets of
$\mathbb{D}$. We begin with the discussion of some properties of families
of analytic functions considered as subsets of ${\mathcal H}(\ID)$.
A univalent function $f$ is called starlike if $f(\ID)$ is a
starlike domain (w.r.t. origin).
Let $\mathcal{S}^*$ denote the class of starlike functions $f\in{\mathcal H}(\ID)$
with $f(0)=0=f '(0)-1$.
Denote by $ {\mathcal K}$ the subclass of functions $f\in{\mathcal H}(\ID)$ with
$f(0)=0=f '(0)-1$ such that $f$ maps $\mathbb{D}$ conformally onto a
convex domain. If
$$\mathcal{B}(M)=\left\{f\in\mathcal{H}(\mathbb{D}):\, f(0)
=f'(0)-1=0 ~\mbox{ and } |zf''(z)|\leq M \right\},
$$
then it is known that $\mathcal{B}(M)\subsetneq \mathcal{K}$ if $0<M\leq 1/2$, and
$\mathcal{B}(M)\subsetneq\mathcal{S}^*$ if $0<M\leq 1$ and the inclusions are sharp.
For a general result we refer to \cite{po-singh02}. In this paper, we are interested
in two subclasses of analytic functions and use the Schwarz lemma as the main tool
in describing the boundary behavior of these two classes of functions.

\begin{nonsec} The Class $\mathcal{B}(\alpha, \beta,M)$ \end{nonsec}\label{subsec1}
Let $\alpha, \beta\in\mathbb{C}$ and $M\in\mathbb{R}$ such that $0<M\leq |\alpha|$
and $|\beta|\leq 1$.
Let $\mathcal{B}(\alpha, \beta,M)$ denote the class of all functions $f$ analytic and univalent
in the unit disk $\ID$,  with $f(0)=0$, $f'(0)=\alpha$, and  $f''(0)=M\beta$
satisfying
\be\label{pvdev7b-eq1}
|zf''(z)|\leq M, \quad z\in\mathbb{D}.
\ee
Note that $\alpha \neq 0$. If $f\in\mathcal{B}(\alpha, \beta,M)$, then we may write
$$zf''(z)=M\omega(z)
$$
for some $\omega \in \mathcal{B}_0$, where
$\mathcal{B}_0$ denotes the class of  functions $\omega$ analytic in
$\mathbb{D}$ such that $|\omega(z)|\leq 1$ in $\mathbb{D}$ and $\omega(0)=0$.
This gives the representation
$$f'(z)-\alpha=f'(z)-f'(0)=M\int_0^1\frac{\omega(tz)}{t}\,dt
$$
so that, by integration,
$$f(z)=\alpha z + Mz\int_0^1\frac{(1-t)\omega(tz)}{t}\,dt.
$$
By the Schwarz lemma, we have $|\omega(z)|\leq |z|$ and so that previous relation
gives that
$$|f'(z)-\alpha|\leq M|z|<M, \quad z\in \ID
$$
which, in particular, shows that functions in $\mathcal{B}(\alpha, \beta,M)$
are univalent in $\mathbb{D}$ if $M\leq |\alpha|$. It is easy to see that functions in
$\mathcal{B}(\alpha,\beta,M)$ are not necessarily univalent if $M>|\alpha|$.
This fact may be demonstrated by, for example, the function
$$f(z)=\alpha z +(M/2)z^2.
$$

Furthermore,  every $f\in \mathcal{B}(\alpha,\beta, M)$  can be associated with a
function $\omega_f$ in $\mathcal{B}_0$ and this association is clearly given by
\be\label{pvdev7b-eq2}
\omega_f(z)=\frac{f''(z)-M\beta}{M-\overline{\beta}f''(z)},\quad z\in\mathbb{D}.
\ee
A simple application of the Schwarz lemma shows that if
$f\in\mathcal{B}(\alpha, \beta, M)$ then one has $|\omega'_f(0)|\leq 1$
which, in particular, gives a restriction on $f'''(0)$. Indeed, it is a simple
exercise to see that
\be\label{pvdev7b-eq2a}
f'''(0)=M(1-|\beta|^2) \omega'_f(0)
\ee
and therefore, with  $\omega'_f(0) = \lambda$, we have
$|f'''(0)|\leq M\left(1-|\beta|^2\right).
$
Using (\ref{pvdev7b-eq2}) and (\ref{pvdev7b-eq2a}), one can obtain by a computation
that
\be\label{pvdev7b-eq02b}
M(1-|\beta|^2)\omega''_f(0)=2M(1-|\beta|^2){\lambda}^2\overline{\beta}+(1+\lambda\overline{\beta})f^{(iv)}(0).
\ee
Also if we let
\be\label{pvdev7a-eq5a}
g(z)=\frac{\frac{\omega_f(z)}{z}-\lambda}
{1-\overline{\lambda}\frac{\omega_f(z)}{z}},\quad\mbox{ for } |\lambda|<1,
\ee
and $g(z)=0$ for $|\lambda|=1$, then $g\in \mathcal{B}_0$, and
we compute that
$$g'(0)=\left\{\ba{ll}
\ds \left.\frac{1}{1-|\lambda|^2}\left(\frac{\omega_f(z)}{z}\right)'\right |_{z=0}
=\frac{1}{1-|\lambda|^2}\left(\frac{\omega''_f(0)}{2}\right)
& \mbox{for $|\lambda|< 1$} \\ [6mm]
\ds 0  & \mbox{for $|\lambda | = 1$.}
 \ea \right .
$$
For convenience, we set $g'(0)=a$.
From (\ref{pvdev7b-eq02b}) we note that for $|\lambda|< 1$,
\be\label{pvdev7a-eq5ab}
|g'(0)|\leq 1
 \Longleftrightarrow
 f^{(iv)}(0)=\frac{2M(1-|\beta|^2)}{1+\lambda\overline{\beta}}
\left[(1-|\lambda|^2)a-{\lambda}^2\overline{\beta}\,\right]
\ee
for some $a\in\overline{\mathbb{D}}$. 
Observe that $a=0$ when $|\lambda |=1$,
and $|a|<1$ if and only if $|\lambda |<1$, according to the Schwarz lemma.

\begin{nonsec} The Class $\mathcal{P}(\alpha,M)$ \end{nonsec}\label{subsec2}
Another class of analytic functions of our interest is defined by
\begin{equation}\label{pvdev7b-eq3a}
\mathcal{P}(\alpha,M)=\left\{f\in\mathcal{H}(\mathbb{D}):\, f(0)=0, f'(0)=\alpha
~\mbox{ and } {\rm Re\,}{zf''(z)}>-M,~z\in\mathbb{D} \right\}
\end{equation}
where $\alpha\in\mathbb{C}\setminus \{0\}$, and $0<M\leq 1/{\log 4}$.
In \cite{ali-samy-singh95}, it has been shown that
$$\mathcal{P}(1,M)\subset\mathcal{S}^* \quad\mbox { for }  0<M\leq \frac{1}{\log 4}.
$$
For any larger value of $M$, functions in $\mathcal{P}(1,M)$ are not necessarily
locally univalent. Later in \cite[Theorem 2.10]{an-par-pon-singh02}
the authors have proved that
$$\mathcal{P}(1,M)\subset\mathcal{S}^*(\alpha)
$$
($0\leq\alpha \leq\frac 12$) whenever
$$0\leq M\leq\frac{1-2\alpha }{2\alpha +\log 4}.
$$   This generalizes the last relation. However, the Herglotz representation
for analytic functions with positive real part
in $\mathbb{D}$ shows that if $f\in\mathcal{P}(\alpha, M) $, then there
exists a unique positive unit measure $\mu$ on $(-\pi, \pi]$ such that
$$\frac{zf''(z)+M}{M}
= \int_{-\pi}^{\pi}\frac{1+ze^{-it}}{1-ze^{-it}}\,d\mu(t),
~\mbox{ i.e. }~
f''(z)=2M \int_{-\pi}^{\pi}\frac{e^{-it}}{1-ze^{-it}}\,d\mu(t).
$$
Integrating from $0$ to $z$ shows that
$$f'(z)=2M \int_{-\pi}^{\pi}
\log\left(\frac{1}{1-ze^{-it}}\right)d\mu(t)+\alpha.
$$
Once again integrating the above from $0$ to $z$ gives the following
representation
$$f(z)=2M\int_{-\pi}^{\pi}\left\{ z-(z-e^{it})\log(1-ze^{-it})\right\}d\mu(t)+\alpha z.
$$
Functions of the above form belong to the class $\mathcal{P}(\alpha, M)$.
Clearly, for every $f\in\mathcal{P}(\alpha, M)$,
there exists an $\omega_f\in\mathcal{B}_0$ such that
\be\label{pvdev7a-eq2}
\omega_f(z)=\frac{zf''(z)}{zf''(z)+2M},\quad z\in\mathbb{D}.
\ee
Note that $|\omega'_f(0)|\leq 1$. By the Schwarz lemma
it is a simple exercise to see that if $f\in\mathcal{P}(\alpha, M)$, then
\be\label{pvdev7a-eq2a}
f''(0)=2M\omega'_f(0)
\ee
and therefore, with  $\omega'_f(0) = \lambda$, we have
$|f''(0)|\leq 2M $. Using (\ref{pvdev7a-eq2}), one can obtain by a
computation that
\be\label{pvdev7a-eq02b}
\omega''_f(0)=\frac{f'''(0)}{M}-2{\lambda}^2.
\ee
Thus, if $g$ is defined by (\ref{pvdev7a-eq5a}) with
$\omega_f(z)$ as in (\ref{pvdev7a-eq2}), then it follows that
\beqq
|g'(0)|\leq 1
& \Longleftrightarrow & f'''(0)=2M((1-|\lambda|^2)a+{\lambda}^2)
\eeqq
for some $a\in\overline{\mathbb{D}}$.
Again we remark that $|g'(0)|=0$ occurs if and only if
$a=\lambda$ with $|\lambda |=1$. 
Also, $|a|<1$ if and only if $|\lambda |<1$, according to the Schwarz lemma.

For $\lambda\in\overline{\mathbb{D}}$ and for each fixed
$z_0\in\mathbb{D}$, we introduce the following sets:
\beqq
\mathcal{S}_1(\lambda)  & = & \Big \{\frac{}{} f\in\mathcal{B}(\alpha, \beta,M)
:\,f'''(0)=M(1-|\beta|^2)\lambda 
\Big \},\\
\mathcal{S}_2(\lambda) & = &  \Big  \{f\in \mathcal{P}(\alpha, M):\, f''(0)=2M\lambda
\Big \},
\eeqq
and
\begin{eqnarray*}
V_j(z_0,\lambda) & = & \{f'(z_0):\, f\in \mathcal{S}_j(\lambda)\} \quad \mbox{for $j=1,2$}.
\end{eqnarray*}
The purpose of the present paper is to determine explicitly the region of variability
$V_j(z_0,\lambda)$ of $ f'(z_0)$ when $f$ ranges over the class $\mathcal{S}_j(\lambda )$
$(j=1,2)$.
Questions of this nature have been discussed recently in
\cite{samy-vasudev1,samy-vasudev3,samy-vasudev-yan3,Yanagihara2}.

\section{The Basic properties of $V_1(z_0,\lambda),$  $V_2(z_0,\lambda)$
and the main results}

For a positive integer $p$, let
$$(\mathcal{S}^*)^p=\{f=f_0^p:\, f_0\in \mathcal{S}^* \}
$$
and recall the following well-known result whose analytic proof is given
in \cite{Yanagihara1} (see also \cite{Goodman's_book,Hallenbeck-Livingston}).

\blem\label{pvdev7b-lem1} Let $f$ be an analytic function in
${\mathbb D}$ with $f(z) = z^p + \cdots $. If
\[ {\rm Re} \,  \left( 1+ z \frac{f''(z)}{f'(z)} \right)
 > 0 , \quad z \in {\mathbb D} ,
\]
then $f \in (\mathcal {S}^*)^p$.
\elem

For the sake of convenience, we use the notation
$V_j(z_0,\lambda)=V_j$ and $\mathcal{S}_j(\lambda )=\mathcal{S}_j$
for $j=1,2$. Now, we begin our investigation by stating certain general
properties of the set $V_j~(j=1,2)$.

\bprop\label{pvdev7b-pro3}
We have
\bee
\item $V_1$ is compact
\item $V_1$ is convex
\item for $|\lambda|=1$ or $z_0=0$,
\be\label{pvdev-7b-eq20}
V_1=\left\{
\frac{M}{\overline{\beta}}\left (z_0-(1-|\beta|^2)
\frac{\log(1+\lambda\overline{\beta}z_0)}{\lambda\overline{\beta}}\right )+\alpha\right\}
\quad\mbox{ if } \beta\neq 0;
\ee
and for $\beta=0$, $V_1$ is obtained as a limiting case from $(\ref{pvdev-7b-eq20})$
\item for $|\lambda|<1$ and $z_0 \in {\mathbb{D}} \setminus \{0 \}$, $V_1$ has
$$
\frac{M}{\overline{\beta}}\left (z_0-(1-|\beta|^2)
\frac{\log(1+\lambda\overline{\beta}z_0)}{\lambda\overline{\beta}}\right )+\alpha
$$
as an interior point.
\eee
\eprop\bpf
(1) First we show that $\mathcal{S}_1$ is compact. For this, we need to
prove that for $\{f_n\}$ in $\mathcal{S}_1$, whenever $f_n\to f$ uniformly on
every compact subset of $\mathbb{D}$, $f\in\mathcal{S}_1$.
We recall that if $f_n\to f$ uniformly on every compact subset of $\mathbb{D}$ then
$f'_n\to f'$ uniformly on every compact subset of $\mathbb{D}$.
Thus, if $f_n\to f$ uniformly,  $f_n(0)\to f(0)$ pointwise which gives $f(0)=0$.
Repeated use of this fact for derivatives, we conclude that
$f'_n(0)\to f'(0)$ pointwise and therefore, $f'(0)=\alpha$.
Similarly, $f''(0)=M\beta$,  $f'''(0)=M(1-|\beta|^2)\lambda$ and
$f^{(iv)}(0)$ satisfies
$$f^{(iv)}(0)=\frac{2M(1-|\beta|^2)}{1+\lambda\overline{\beta}}\left[(1-|\lambda|^2)a
-{\lambda}^2\overline{\beta}\,\right].
$$
Also, we have
$$|zf_n''(z)|\to |zf''(z)|
$$
uniformly on every compact subset of $\mathbb{D}$.
Since each member of the sequence $\{f_n\}$ is in $\mathcal{S}_1$,
it follows that  $|zf''(z)|\leq M$. We conclude that
$\mathcal{S}_1$ is compact.

Finally, for fixed $z_0\in\mathbb{D}$, define
$\psi :\, \mathcal{S}_1\to V_1$  by
$$
\psi(f)(z_0)=f'(z_0).
$$
Clearly $\psi$ is continuous. Thus, $V_1$ is compact.


(2) If $f_0$ and $f_1$ belong to $\mathcal{S}_1$, then, for $0\leq t \leq 1$, the
function $f_t$ defined by
$$f_t(z)=\int_0^z \{(1-t)f'_0(\zeta)+tf'_1(\zeta)\}\,d\zeta
$$
also belongs to $\mathcal{S}_1$. Clearly,
we have
$$f'_t(z)=(1-t)f'_0(z)+tf'_1(z),
$$
and the convexity of $V_1$ is evident.

(3) If $z_0=0$, (\ref{pvdev-7b-eq20}) trivially holds. If $|\lambda|=1$,
then from (\ref{pvdev7b-eq2a}) we see that $|\omega '_f(0)|=1$
and therefore
it follows from the Schwarz lemma that $\omega_f(z)=\lambda z$, which
by (\ref{pvdev7b-eq2}) gives that
$$ \frac{zf''(z)}{M}=\frac{(\lambda z+\beta)z}{1+\overline{\beta}\lambda z},
\quad\mbox{ for }~ |\beta|\leq 1.
$$
By integrating the above from $0$ to $z_0$ we see that
\begin{align*}
f'(z_0)&=M\int_0^{z_0}\frac{\lambda\zeta+\beta}{1+\overline{\beta}\lambda\zeta}
\,d\zeta +\alpha \quad\mbox{if }~ \beta\neq 0\\
&=\frac{M}{\overline{\beta}}\int_0^{z_0}
\left[\left( 1-\frac{1}{1+\overline{\beta}\lambda\zeta}\right)
+\frac{\beta}{\lambda}\left( \frac{\lambda\overline{\beta}}
{1+\overline{\beta}\lambda\zeta}\right)\right]d\zeta+\alpha \quad\mbox{if }~ \beta\neq 0
\end{align*}
and a computation gives
$$f'(z_0) =\left\{\ba{ll}
\ds\frac{Mz_0}{\overline{\beta}}-\frac{M}{\lambda\overline{\beta}^2}
(1-|\beta|^2)\log(1+\lambda\overline{\beta}z_0)+\alpha
& \mbox{if $\beta\neq 0$} \\ [6mm]
\ds\alpha+\frac{M\lambda}{2}z_0^2
 & \mbox{if $\beta=0$}
 \ea \right .
$$
and thus,
$$V_1=\left\{\frac{Mz_0}{\overline{\beta}}-\frac{M}{\lambda\overline{\beta}^2}
(1-|\beta|^2)\log(1+\lambda\overline{\beta}z_0)+\alpha \right\}\quad\mbox{if }
~\beta\neq 0.
$$
We remark that
$$
\lim_{\beta\to 0}\left\{\frac{Mz_0}{\overline{\beta}}-\frac{M}{\lambda\overline{\beta}^2}
(1-|\beta|^2)\log(1+\lambda\overline{\beta}z_0)+\alpha \right\}
=\alpha+\frac{M\lambda}{2}z_0^2
$$
and, therefore for $\beta=0$,
$$V_1=\left\{\alpha+\frac{M\lambda}{2}z_0^2 \right\}.
$$
Hence, the extremal function in $\mathcal{S}_1$ for $|\lambda|=1$ is of the
form
$$f(z) =\left\{\ba{ll}
\ds
\frac{M}{\overline{\beta}}\left\{\frac{z^2}{2}-\frac{(1-|\beta|^2)}{\lambda\overline{\beta}}
\left[\left(z+\frac{1}{\lambda\overline{\beta}}\right)
\log(1+\lambda\overline{\beta}z)-z\right]\right\}+\alpha z
& \\ [6mm]
\qquad \qquad \qquad \mbox{ if $0<|\beta|\leq 1$} & \\
\ds \alpha z+\frac{M\lambda}{6}z^3\quad \mbox{if $\beta=0$}.
 &
 \ea \right .
$$

(4) For $\lambda\in\mathbb{D}$, we
let
\begin{eqnarray}\label{pvdev7b-delta}
\delta(z,\lambda)
& = & \frac{z+\lambda}{1+\overline{\lambda}z}.
\end{eqnarray}
A simplification of (\ref{pvdev7a-eq5a}) with $g(z)=az$ ($|a|<1$) leads to
\begin{eqnarray}\label{pvdev7b-extremal}
F_{a,\lambda}(z):=f(z)
 & = &
 \int_0^z\left\{\int_0^{\zeta_{2}}
\left(\frac{M(\delta(a\zeta_1,\lambda)\zeta_1+\beta)}
{1+\overline{\beta}\delta(a\zeta_1,\lambda)\zeta_1}\right)d\zeta_1\right\}d\zeta_2
+\alpha z,
\end{eqnarray}
where $z\in\mathbb{D}$.

First we claim that $F_{a,\lambda}\in \mathcal{S}_1$.
In fact, with the aid of (\ref{pvdev7b-extremal}) we easily get
$$zF''_{a,\lambda}(z)
=\frac{M(\delta (az,\lambda)z+\beta)z}{1+\overline{\beta}\delta(az,\lambda)z}.
$$
As $\delta (az,\lambda)$ lies in the unit disk $\mathbb{D}$,
$F_{a,\lambda}\in\mathcal{S}_1$. Further, one may verify that
\be\label{pvdev7b-eq3}
\omega_{F_{a,\lambda}}(z)=z\delta(az,\lambda).
\ee
Next our claim is that for a fixed $\lambda\in\mathbb{D}$ and
$z_0\in\mathbb{D}\setminus\{0\}$,
\begin{eqnarray*}
\mathbb{D}\ni a \mapsto  F'_{a,\lambda}(z_0)=
\int_0^{z_0} \left(\frac{M(\delta(a\zeta,\lambda)\zeta+\beta)}
{1+\overline{\beta}\delta(a\zeta,\lambda)\zeta}\right)d\zeta+\alpha,
\end{eqnarray*}
is a non-constant analytic function of $a\in\mathbb{D}$, and hence is an open mapping.

Finally, we claim that the mapping
${\mathbb D} \ni a\mapsto  F'_{a,\lambda}(z_0) $
is a non-constant analytic function of
$a$ for each fixed $z_0 \in {\mathbb D} \backslash \{ 0 \}$ and
$\lambda\in\mathbb{D}$. For this, we put
$$ h_1(z) = \left .\frac{3}{M(1-|\beta|^2)(1-|\lambda|^2)}
\frac{\partial}{\partial a}\left\{\frac{}{} F'_{a,\lambda} (z)\right\}\right |_{a=0}
$$
and obtain that
$$h_1(z)  =  3\int_0^z \frac{{\zeta}^2}{(1+\overline{\beta}\lambda \zeta)^2}\,d\zeta
=z^3+\cdots
$$
which gives
$$ {\rm Re} \,\left\{\frac{zh_1''(z)}{h_1'(z)}\right\} = 2\,
{\rm Re} \,\left\{\frac{1}{1+\overline{\beta}\lambda z}\right\}\geq
\frac{2}{1+|\lambda |\, |\beta|}> 1, \quad z\in \ID .
$$
By Lemma \ref{pvdev7b-lem1}, there exists a function $h_0\in
\mathcal{S}^*$ with $h_1=h_0^3$. The univalence of $h_0$ and
$h_0(0)=0$ imply  that $h_1(z_0)\neq 0$ for $z_0 \in {\mathbb
D}\setminus \{0\}$. Consequently, the mapping
${\mathbb D} \ni a\mapsto F'_{a,\lambda}(z_0) $ is a non-constant analytic  function
of $a$. Thus
$$F'_{0,\lambda}(z_0)=\int_0^{z_0}M\left(\frac{\lambda\zeta+\beta}
{1+\overline{\beta}\lambda\zeta}\right)d\zeta +\alpha
$$
is an interior point of
$\{ F'_{a,\lambda}(z_0):\,a\in\mathbb{D}\}\subset V_1$.
\epf



We have the following analog of Proposition \ref{pvdev7b-pro3}
for functions in $\mathcal{P}(\alpha,M)$.

\bprop\label{pvdev7a-pro3}
We have
\bee
\item $V_2$ is compact
\item $V_2$ is convex
\item for $|\lambda|=1$ or $z_0=0$,
\be\label{pvdev-7a-eq20}
V_2=\{-2M\log(1-\lambda z_0)+\alpha\}
\ee
\item for $|\lambda|<1$ and $z_0 \in {\mathbb{D}} \setminus \{0 \}$, $V_2$ has $-2M\log(1-\lambda z_0)+\alpha$
as an interior point
\eee
\eprop\bpf
Part $(1)$ and $(2)$ follows exactly as in the proof of Proposition \ref{pvdev7b-pro3} and
so we omit the details.

%
%
%

(3) If $z_0=0$, (\ref{pvdev-7a-eq20}) trivially holds. If $|\lambda|=1$, then
for $f\in\mathcal{P}(\alpha,M)$ it follows that $f''(0)=2M\omega'_f(0)$ and
from (\ref{pvdev7a-eq2a}) we see that $|\omega'_f(0)|=1$ and therefore it follows from the Schwarz lemma
that $\omega_f(z)=\lambda z$, which by (\ref{pvdev7a-eq2}) gives
$$ zf''(z)+M=\frac{M(1+\lambda z)}{1-\lambda z}\quad\mbox{ or }
f''(z)=\frac{2M\lambda}{1-\lambda z}.
$$
By integrating from $0$ to $z_0$ we see that
$$f'(z_0)=-2M\log(1-\lambda z_0)+\alpha.
$$
Thus
$$ V_2=\left\{-2M\log(1-\lambda z_0)+\alpha\right\}.
$$
Hence, the extremal function in $\mathcal{S}_2$ for $|\lambda|=1$ is of the form
$$
f(z)=2M\left(\frac{1}{\lambda}-z\right)\log\left(\frac{1}{1-\lambda z}\right)+(\alpha-2M)z.
$$

(4) Let $\lambda\in\mathbb{D}$ and $a\in\mathbb{D}$.
A simple computation  as before helps to introduce
\begin{eqnarray}\label{pvdev7a_extremal}
H_{a,\lambda}(z):=f(z)
 & = &
 \int_0^z\left\{\int_0^{\zeta_{2}}
\frac{2M\delta(a\zeta_1,\lambda)}
{1-\delta(a\zeta_1,\lambda)\zeta_1}\,d\zeta_1\right\}d\zeta_2+\alpha z, \quad z\in\mathbb{D},
\end{eqnarray}
where $\delta(z,\lambda)$ is defined by (\ref{pvdev7b-delta}).
From this we see that $H_{a,\lambda}\in \mathcal{S}_2$
and
\be\label{pvdev7a-eq3}
\omega_{H_{a,\lambda}}(z)=z\delta(az,\lambda).
\ee
For a fixed $\lambda\in\mathbb{D}$ and $z_0\in\mathbb{D}\setminus\{0\}$,
the function
\begin{eqnarray*}\\
\mathbb{D}\ni a \mapsto  H'_{a,\lambda}(z_0)=
\int_0^{z_0} \frac{2M\delta(a\zeta,\lambda)}
{1-\delta(a\zeta,\lambda)\zeta}\,d\zeta+\alpha
\end{eqnarray*}
is a non-constant analytic function of $a\in\mathbb{D}$, and hence is an open mapping.
We claim that the mapping ${\mathbb D} \ni a\mapsto  H'_{a,\lambda}(z_0)$
is a non-constant analytic function of $a$ for each fixed
$z_0 \in {\mathbb D} \backslash \{ 0 \}$ and $\lambda\in\mathbb{D}$. For this we let
\begin{eqnarray*}
h_2(z) & = & \left .\frac{1}{M(1-|\lambda|^2)}
\frac{\partial}{\partial a}\left\{\frac{}{} H'_{a,\lambda} (z)\right\}\right |_{a=0}
\end{eqnarray*}
so that
$$h_2(z) = 2\int_0^z \frac{\zeta}{(1-\lambda\zeta)^2}\,d\zeta =z^2+\cdots .
$$
This gives
$$ {\rm Re} \,\left\{\frac{zh_2''(z)}{h_2'(z)}\right\} =
{\rm Re} \,\left\{\frac{1+\lambda z}{1-\lambda z}\right\}>0, \quad z\in \ID .
$$
By Lemma \ref{pvdev7b-lem1} there exists a function $h_0\in{\mathcal S}^*$ with $h_2=h_0^2$
so that the univalence of $h_0$ and $h_0(0)=0$ imply  that $h_2(z_0)\neq 0$ for
$z_0 \in \mathbb {D}\setminus \{0\}$. Consequently, the mapping ${\mathbb D} \ni
a\mapsto H'_{a,\lambda}(z_0) $ is a non-constant analytic  function
of $a$. Thus,
$$H'_{0,\lambda}(z_0)= -2M\log(1-\lambda z_0)+\alpha
$$
is an interior point of
$\{H'_{a,\lambda}(z_0):\,a\in\mathbb{D}\}\subset V_2$.

\epf


For each $j=1,2$,  $V_j$ is a compact convex subset of $\mathbb{C}$
and has nonempty interior and therefore,  the boundary
$\partial{V_j}$ is a Jordan curve and $V_j$ is the union of $\partial{V_j}$
 and its inner domain. We now state our main results.

\bthm\label{pvdev7b-th1}
For $\lambda\in\mathbb{D}$, the boundary
$\partial{V_1}$ is the Jordan curve given by
\begin{eqnarray*}
(-\pi,\pi]\ni \theta \mapsto F'_{e^{i\theta},\lambda}(z_0)
 & = &
\int_0^{z_0}
\left(\frac{M(\delta(e^{i\theta}\zeta,\lambda)\zeta+\beta)}
{1+\overline{\beta}\delta(e^{i\theta}\zeta,\lambda)\zeta}\right)d\zeta+\alpha.
\end{eqnarray*}
Here $\alpha, \beta\in\mathbb{C}$ and $M\in\mathbb{R}$  with  $0<M\leq |\alpha|$
and $|\beta|\leq 1$.
If $f'(z_0)=F'_{e^{i\theta},\lambda}(z_0)$ for some $f\in\mathcal{S}_1$,
$z_0\in\mathbb{D}\setminus \{0\}$,
and $\theta\in (-\pi,\pi]$, then $f(z)= F_{e^{i\theta},\lambda}(z)$.
\ethm

\bthm\label{pvdev7a-th2} For $\lambda\in\mathbb{D}$,
the boundary $\partial{V_2}$ is the Jordan curve given by
\begin{eqnarray*}
(-\pi,\pi]\ni \theta \mapsto H'_{e^{i\theta},\lambda}(z_0)
 & = & \int_0^{z_0}
\frac{2M\delta(e^{i\theta}\zeta,\lambda)}
{1-\delta(e^{i\theta}\zeta,\lambda)\zeta}\,d\zeta+\alpha
\end{eqnarray*}
where $\alpha\in\mathbb{C}\setminus \{0\}$ and $0<M\leq 1/{\log 4}$.
If $f'(z_0)=H'_{e^{i\theta},\lambda}(z_0)$ for some $f\in \mathcal{S}_2$,
$z_0\in\mathbb{D}\setminus \{0\}$
and $\theta\in (-\pi,\pi]$, then $f(z)= H_{e^{i\theta},\lambda}(z)$.
\ethm
\section{Preparation for the proof of Theorem \ref{pvdev7b-th1}}
We exclude the case $|\lambda|=1$ from the following result
as this follows from our earlier discussion.

\bprop\label{pvdev7b-pro1}
For $f\in\mathcal{S}_1(\lambda )$ with
$|\lambda |<1$,
 we have
\be\label{pvdev7b-eq4}
\left|f''(z)-c_1(z,\lambda)\right|\leq r_1(z,\lambda), \quad z\in\mathbb{D},
\ee
where
\begin{eqnarray*}
c_1(z,\lambda) & = &
\frac{M(1-|z|^2)\left\{\beta(1+|z|^2)+{\beta}^2\overline{\lambda}\overline{z}+\lambda
z\right\}}
{\left(1-|\beta|^2|z|^4\right)-\left(1-|\beta|^2\right)|\lambda|^2|z|^2
+2\left(1-|z|^2\right){\rm Re\,}(\overline{\beta}\lambda z)},\quad \mbox{and}\\
r_1(z,\lambda) & = &
\frac{(1-|\lambda|^2)\left(1-|\beta|^2\right)|z|^2}
{\left(1-|\beta|^2|z|^4\right)-\left(1-|\beta|^2\right)|\lambda|^2|z|^2
+2\left(1-|z|^2\right){\rm Re\,}(\overline{\beta}\lambda z)}.
\end{eqnarray*}
For each $z\in\mathbb{D}\setminus\{0\}$, equality holds if and only if
$f=F_{e^{i\theta},\lambda}$ for some $\theta\in\mathbb{R}$.
\eprop
\bpf
Let $f\in\mathcal{S}_1$. Then (\ref{pvdev7b-eq2}) holds with $\omega_f(\in\mathcal{{B}}_0)$
and  $\omega'_f(0)=\lambda$. It follows from the Schwarz lemma that
\be\label{pvdev7b-eq5}
\left|\frac{\frac{\omega_f(z)}{z}-\lambda}{1-{\overline\lambda}\frac{\omega_f(z)}{z}}\right|\leq |z|.
\ee
From  (\ref{pvdev7b-eq1}) and (\ref{pvdev7b-eq2}) this is easily seen to be equivalent to
\be\label{pvdev7b-eq6}
\left|\frac{f''(z)-A_1(z,\lambda)}{f''(z)+B_1(z,\lambda)}
\right| \leq |z|\,|\tau_1(z,\lambda)|,
\ee
where
\be\label{pvdev7b-eq7}
\left \{ \ba{l}
\ds A_1(z,\lambda) = \frac{M(\beta+\lambda z)}{1+\overline{\beta}\lambda z},\\
\ds B_1(z,\lambda) =
-\frac{M(z+\overline{\lambda}\beta)}{\overline{\beta}z+\overline{\lambda}}, \\
\ds \tau_1(z,\lambda)
 = \frac{\overline{\beta}z+\overline{\lambda}}{1+\overline{\beta}\lambda z}.
\ea
\right .
\ee
Further, a computation shows that the inequality (\ref{pvdev7b-eq6}) is equivalent to
\be\label{pvdev7b-eq8}
\left|f''(z)-\frac{A_1(z,\lambda)+|z|^2|\tau_1(z,\lambda)|^2 B_1(z,\lambda)}
{1-|z|^2|\tau_1(z,\lambda)|^2}\right| \leq
\frac{|z|\,|\tau_1(z,\lambda)|\,|A_1(z,\lambda)+B_1(z,\lambda)|}{1-|z|^2|\tau_1(z,\lambda)|^2}.
\ee
Now we have
\begin{align}
1 - |z|^2|\tau_1(z,\lambda)|^2
 & =1-|z|^2\left|\frac{\overline{\beta}z+\overline{\lambda}}{1+\overline{\beta}\lambda z}\right|^2\nonumber\\
 & = \frac{{\left(1-|\beta|^2|z|^4\right)-\left(1-|\beta|^2\right)|\lambda|^2|z|^2
+2\left(1-|z|^2\right){\rm Re\,}(\overline{\beta}\lambda z)}}
{\left|1+\overline{\beta}\lambda z\right|^2}\nonumber,
\end{align}
\begin{align*}
A_1(z,\lambda)  +  B_1(z,\lambda)
& = \frac{M(\beta+ \lambda z)}{1+\overline{\beta}\lambda z}
-\frac{M(z+\overline{\lambda}\beta)}{\overline{\beta}z+\overline{\lambda}}\\\nonumber
& =\frac{M(1-|\lambda|^2)(1-|\beta|^2)z}{\left(1+\overline{\beta}\lambda z\right)
\left(\overline{\beta}z+\overline{\lambda}\right)}\nonumber
\end{align*}
and
\begin{align*}
 A_1(z,\lambda)&+|z|^2|\tau_1(z,\lambda)|^2 B_1(z,\lambda)\\
& =\frac{M(\beta+\lambda z)}{1+\overline{\beta}\lambda z} -
|z|^2\left|\frac{\overline{\beta}z+\overline{\lambda}}
{1+\overline{\beta}\lambda z}\right|^2
\left(\frac{M(z+\overline{\lambda}\beta)}
{\overline{\beta}+\overline{\lambda}}\right)\\
& =\frac{M(1-|z|^2)\left\{\beta(1+|z|^2)+{\beta}^2\overline{\lambda}\overline{z}+\lambda z\right\}}
{\left|1+\overline{\beta}\lambda z\right|^2}\nonumber.
\end{align*}
An easy calculation yields that
$$
\frac{A_1(z,\lambda)+|z|^2|\tau_1(z,\lambda)|^2 B_1(z,\lambda)}
{1-|z|^2|\tau_1(z,\lambda)|^2}= c_1(z,\lambda)
$$
and
$$
\frac{|z|\,|\tau_1(z,\lambda)|\,|A_1(z,\lambda)+B_1(z,\lambda)|}{1-|z|^2|\tau_1(z,\lambda)|^2}
= r_1(z,\lambda).
$$
Now the inequality (\ref {pvdev7b-eq4}) follows from these equalities and
(\ref{pvdev7b-eq8}).

It is easy to see that the equality occurs for $z\in\mathbb{D}$ in (\ref {pvdev7b-eq4})
if and only if equality occurs in (\ref{pvdev7b-eq5}). Thus
the equality in (\ref {pvdev7b-eq4}) holds whenever
$f=F_{e^{i\theta},\lambda}$ for some $\theta\in\mathbb{R}$.
Conversely if the equality occurs for some
$z\in\mathbb{D}\setminus\{0\}$ in (\ref {pvdev7b-eq4}), then the
equality must hold in (\ref{pvdev7b-eq8}) and hence (\ref{pvdev7b-eq5}) holds.
Thus, from the Schwarz lemma, there  exists a $\theta\in\mathbb{R}$ such that
$\omega_f(z)=z\delta(e^{i\theta}z,\lambda)$ for all
$z\in\mathbb{D}$. This implies $f=F_{e^{i\theta},\lambda}$.
\epf

The case  $\lambda =0$ of Proposition \ref{pvdev7b-pro1}
gives the following  result.

\bcor\label{pvdev7b-cor3}
Let $f\in\mathcal{S}_1(0)$. Then we have
$$\left|f''(z) -\frac{M\beta(1-|z|^4)}{1-|\beta|^2|z|^4}\right |
\leq \frac{\left(1-|\beta|^2\right)|z|^2}{1-|\beta|^2|z|^4},\quad z\in\mathbb{D}.
$$
For each $z\in\mathbb{D}\setminus\{0\}$, equality holds if and only if
$f=F_{e^{i\theta},0}$ for some $\theta\in\mathbb{R}$. Here $F_{e^{i\theta},0}$
is defined in Theorem \ref{pvdev7b-th1}.
\ecor
For $|\beta|=1$, by Corollary \ref{pvdev7b-cor3}, functions in
$\mathcal{S}_1(0)$ must
satisfy
$$|f''(z)-M\beta|=0
$$
which gives
$$f(z)=\alpha z+\frac{M\beta}{2}z^2.
$$

\bcor\label{pvdev7b-cor1}
Let $\gamma:\,z(t)$, $0\leq t\leq 1$, be a $C^1$-curve in
$\mathbb{D}$ with $z(0)=0$ and $z(1)=z_0$. Then we have
$$
V_1\subset\overline{\mathbb{D}}(C_1(\lambda, \gamma),R_1(\lambda, \gamma))
=\left\{w\in\mathbb{C}:\, |w-C_1(\lambda, \gamma)|\leq R_1(\lambda, \gamma)\right\},
$$
where
$$
C_1(\lambda, \gamma)=\alpha + \int_0^1 c_1(z(t),\lambda)z'(t)\,dt
~\mbox{ and }~
R_1(\lambda, \gamma)=\int_0^1 r_1(z(t),\lambda)|z'(t)|\,dt.
$$
\ecor \bpf
Since for $f\in\mathcal{S}_1$ we have
$$\int_0^1 f''(z(t))z'(t)\,dt= f'(z_0)-f'(0)=f'(z_0)-\alpha,
$$
it follows from Proposition \ref{pvdev7b-pro1} that
\begin{eqnarray*}
|f'(z_0)-C_1(\lambda, \gamma)|
& = & \left| f'(z_0)-\alpha -\int_0^1 c_1(z(t),\lambda)z'(t)\,dt\right|\\
& = & \left|\int_0^1\left\{f''(z(t))-c_1(z(t),\lambda)\right\}z'(t)\,dt\right|\\
&\leq& \int_0^1 r_1(z(t),\lambda)|z'(t)|\,dt = R_1(\lambda, \gamma).
\end{eqnarray*}
As $f'(z_0)\in V_1$ was arbitrary, the conclusion follows.
\epf
\blem\label{pvdev7b-lem2}
For $\theta\in\mathbb{R}$, $\lambda\in\mathbb{D}$ and $\beta\in\overline{\mathbb{D}}$,
the function
$$G(z)= \int_0^z \frac{e^{i\theta}{\zeta}^2}{\left\{1+\left(\overline{\lambda}e^{i\theta}
+\overline{\beta}\lambda\right)\zeta+\overline{\beta}e^{i\theta}{\zeta}^2\right\}^2}\,d\zeta,
\quad z\in\mathbb{D},
$$
has a zero of order three at the origin and no zeros elsewhere in $\mathbb{D}$.
Furthermore, there exists a starlike univalent function $G_0$ in $\mathbb{D}$
such that $G=3^{-1}e^{i\theta}G^3_0$ and $G_0(0)= G'_0(0)-1=0$.
\elem
\bpf
We first prove that
\begin{equation}\label{pvdev7b-eq22}
{\rm Re} \left\{\frac{z G''(z)}{G'(z)}\right\}>-1,\quad z\in\mathbb{D}.
\end{equation}

If $\beta=0$, then a simple computation gives (\ref{pvdev7b-eq22}).

If $0<|\beta|\leq 1$, then it is easy to see that
$$1+(\overline\lambda e^{i\theta}+\overline{\beta}\lambda)z+
\overline{\beta}e^{i\theta}z^2
= \left(1-\frac{z}{z_1}\right)\left(1-\frac{z}{z_2}\right),
$$
where
$$z_1z_2=\frac{e^{-i\theta}}{\overline{\beta}} ~\mbox{ and }~
z_1+z_2 = \frac{\overline{\lambda}
+ e^{-i\theta}\overline{\beta}\lambda}{\overline{\beta}}.
$$
We note that $|z_1|\,|z_2|=1/{|\beta|}$. We need to show that none of $z_1$
and $z_2$ lie in the punctured unit disk $\ID \backslash \{0\}$.

Suppose first that $0<|\beta |<1$. Then either
both $z_1$ and $z_2$ lie outside the unit circle,  or else one lies inside
while the other lies outside the unit circle.
We claim that the later case cannot occur. On the
contrary, without loss of generality, we may assume that
$$ |z_1|<1,~ |z_2|>1.
$$
Also let $z_1=re^{i\phi}$, for some $r<1$ and $\phi\in\mathbb{R}$. Then
$$z_2=\frac{e^{-i(\theta+\phi)}}{r\overline{\beta}}
$$
and so the expression
$$z_1+z_2=\frac{\overline{\lambda}+ e^{-i\theta}\overline{\beta}\lambda}{\overline{\beta}}
$$
simplifies to an equivalent form
\begin{equation}\label{pvdev7b-eq21}
\zeta +\frac{1}{|\beta|\zeta}=-\left(\omega+\frac{1}{|\beta|\omega}\right)
\end{equation}
with $|\zeta|=1$ and $|\omega|=|\lambda|$. We may rewrite (\ref{pvdev7b-eq21}) as
$$ (\zeta+\omega)\left(1+\frac{1}{|\beta|\omega\zeta}\right)=0
$$
which is a contradiction, because
this equation has no solution when $|\zeta|=1$ and $|\omega|=|\lambda|$.
 We conclude that $|z_1|>1$ and $|z_2|>1$.

If $|\beta|=1$, then $|z_1|\,|z_2|=1$ so that
either $|z_1|=1$ and $|z_2|=1$, or
$|z_1|<1$ and $|z_2|>1$, or $|z_1|>1$ and $|z_2|<1$ holds.
Again we see that the last two cases cannot occur. Indeed,
on the contrary, we may (without loss of generality) assume that
$$|z_1|<1,~ |z_2|>1.
$$
Then the expression for $z_1+z_2$ simplifies to the form
\begin{equation}\label{pvdev7b-eq21c}
\zeta +\frac{1}{\zeta}=-{\rm Re \,}(\overline{\lambda}e^{i\psi}),
\end{equation}
with  $|\zeta|=r<1$ and $\psi\in\mathbb{R}$. Now the set of complex
numbers described by the right
hand side of (\ref{pvdev7b-eq21c}) forms a subset of real numbers lying in the
line segment $(-1, 1)$ whereas the set of  complex numbers described
by the left hand side of (\ref{pvdev7b-eq21c}) lies out side of the ellipse
$$ \frac{u^2}{(1/4)(r+1/r)^2}+\frac{v^2}{(1/4)(r-1/r)^2}=1.
$$
Since the above two sets of complex numbers are disjoint, we arrive at
a contradiction. Hence, we conclude that $|z_1|=|z_2|=1$.

Finally, as $|z_1|\geq 1$ and $|z_2|\geq 1$,  a simple calculation shows that
$${\rm Re}\,\left\{\frac{zG''(z)}{G'(z)}\right\}
= {\rm Re}\, \left(\frac{1+z/z_1}{1-z/z_1}\right)
+ {\rm Re}\, \left(\frac{1+z/z_2}{1-z/z_2}\right)>0, \quad z\in\mathbb{D}.
$$

Applying Lemma \ref{pvdev7b-lem1} to $3e^{-i \theta}G(z)$ with
$p=3$ there exists a $G_0\in\mathcal{S}^*$ such that $G=3^{-1}e^{i\theta}G_0^3$.
This completes the proof.
\epf

\bprop\label{pvdev7b-pro2}
Let $z_0\in\mathbb{D}\setminus \{0\}$. Then for $\theta\in(-\pi,\pi]$
we have $F'_{e^{i\theta},\lambda}(z_0)\in\partial{V_1}$. Furthermore if
$f'(z_0)= F'_{e^{i\theta},\lambda}(z_0)$ for some
$f\in\mathcal{S}_1$ and $\theta\in(-\pi,\pi]$, then $f= F_{e^{i\theta},\lambda}$.
\eprop\bpf
From (\ref{pvdev7b-extremal}) we easily obtain that
$$F''_{a, \lambda}(z)
=\frac{M(\delta (az,\lambda)z+\beta)}{1+\overline{\beta}\delta(az,\lambda)z}
= \frac{M\left\{(az+\lambda)z+\beta(1+\overline{\lambda}az)\right\}}
{1+\left(\overline{\lambda}a+\overline{\beta}\lambda\right)z+\overline{\beta}az^2}.
$$
Thus we have from (\ref{pvdev7b-eq7})
$$F''_{a, \lambda}(z)-A_1(z,\lambda)
=\frac{M(1-|\lambda|^2)\left(1-|\beta|^2\right)az^2}
{\left(1+\left(\overline{\lambda}a+\overline{\beta}\lambda\right)z+\overline{\beta}az^2\right)
\left(1+\overline{\beta}\lambda z\right)},
$$
$$
F''_{a, \lambda}(z)+B_1(z,\lambda)
=\frac{-M(1-|\lambda|^2)(1-|\beta|^2)z}
{\left(1+\left(\overline{\lambda}a+\overline{\beta}\lambda\right)z
+\overline{\beta}az^2\right)\left(\overline{\beta}z+\overline{\lambda}\right)}
$$
and hence

\vspace{8pt}
\noindent
$F''_{a, \lambda}(z)-c_1(z,\lambda)$
\begin{eqnarray*}
 & = & F''_{a, \lambda}(z)-\frac{A_1(z,\lambda)+
|z|^2|\tau_1(z,\lambda)|^2 B_1(z,\lambda)}
{1-|z|^2|\tau_1(z,\lambda)|^2}\\
 &  =  &
\frac{1}{1-|z|^2 |\tau_1(z,\lambda)|^2}\left\{(F''_{a, \lambda}(z)
-A_1(z,\lambda))-|z|^2|\tau_1(z,\lambda)|^2
\left(F''_{a, \lambda}(z)+B_1(z,\lambda)\right)\right\}\\
 & = & \frac{M(1-|\lambda|^2)(1-|\beta|^2)az^2 }
{\left\{(1-|\beta|^2|z|^4)-\left(1-|\beta|^2\right)|\lambda|^2|z|^2
+2(1-|z|^2){\rm Re\,}(\overline{\beta}\lambda z)\right\}}
\frac{\overline{K(a,z)}}{K(a,z)},
\end{eqnarray*}
where
$$ K(a,z)= 1+\left(\overline{\lambda}a+\overline{\beta}\lambda\right)z+\overline{\beta}az^2.
$$
Substituting $a=e^{i\theta}$, we find that

\vspace{8pt}
\noindent
$F''_{e^{i\theta}, \lambda}(z)-c_1(z,\lambda)$
\begin{eqnarray*}
 & = &
\frac{M(1-|\lambda|^2)(1-|\beta|^2)e^{i\theta}z^2}
{\left\{(1-|\beta|^2|z|^4)-\left(1-|\beta|^2\right)|\lambda|^2|z|^2
+2(1-|z|^2){\rm Re\,}(\overline{\beta}\lambda z)\right\}}
\frac{\overline{K(e^{i\theta},z)}}{K(e^{i\theta},z)}
\\
  & = &
\frac{M(1-|\lambda|^2)(1-|\beta|^2)e^{i\theta}z^2}
{\left\{\left(1-|\beta|^2|z|^4\right)-\left(1-|\beta|^2\right)|\lambda|^2|z|^2
+2\left(1-|z|^2\right){\rm Re\,}(\overline{\beta}\lambda z)\right\}}
\frac{\left|K(e^{i\theta},z)\right|^2}{(K(e^{i\theta},z))^2}
\\
 & = &
r_1(z,\lambda)\frac{e^{i\theta}z^2}{|z|^2}
\frac{\left|1+\left(\overline{\lambda}e^{i\theta}+\overline{\beta}\lambda\right)z
+\overline{\beta}e^{i\theta}z^2\right|^2}{\left\{1+\left(\overline{\lambda}e^{i\theta}
+\overline{\beta}\lambda\right)z+\overline{\beta}e^{i\theta}z^2\right\}^2}.
\end{eqnarray*}
From Lemma \ref{pvdev7b-lem2}, we may rewrite the last expression as
\be\label{pvdev7b-eq9}
F''_{e^{i\theta}, \lambda}(z)-c_1(z,\lambda)
= r_1(z,\lambda)\frac{G'(z)}{|G'(z)|},
\ee
where $G(z)$ is defined as in Lemma \ref{pvdev7b-lem2}. According to
Lemma \ref{pvdev7b-lem2}, the function $G_0$ defined by
$ G=3^{-1}e^{i\theta}G^3_0$ is starlike. As a consequence,
for any $z_0\in\mathbb{D}\setminus\{0\}$, the line
segment joining $0$ and $G_0(z_0)$ entirely lies in $G_0(\mathbb{D})$.
Introduce $\gamma_0$ by
\be\label{pvdev7b-eq10}
\gamma_0: \,z(t)=G_0^{-1}(tG_0(z_0)), \quad 0\leq t \leq 1.
\ee
From the representation of $G$, we obtain
$$G(z(t))= 3^{-1}e^{i\theta}G_0(z(t))^3= 3^{-1}e^{i\theta}(tG_0(z_0))^3
=t^3 G(z_0)
$$
and so,
\be\label{pvdev7b-eq11}
G'(z(t))z'(t)=3t^2G(z_0), \quad t\in [0,1].
\ee
Using this and (\ref{pvdev7b-eq9}) we deduce that
\begin{eqnarray}
\label{pvdev7b-eq12}
 F'_{e^{i\theta},\lambda}(z_0)-C_1(\lambda, \gamma_0)
 & = & \int_0^1\left\{F''_{e^{i\theta}, \lambda}(z(t))-c_1(z(t),\lambda)\right\}z'(t)\,dt\\
 & = &  \int_0^1 r_1(z(t),\lambda)\frac{G'(z(t))z'(t)}{|G'(z(t))z'(t)|}|z'(t)|\,dt\nonumber \\
 & = & \frac{G(z_0)}{|G(z_0)|}\int_0^1 r_1(z(t),\lambda)|z'(t)|\,dt\nonumber \\
 & = & \frac{G(z_0)}{|G(z_0)|}R_1(\lambda, \gamma_0)\nonumber
\end{eqnarray}
which means that $F'_{e^{i\theta},\lambda}(z_0)\in\partial{\overline{\mathbb{D}}}
(C_1(\lambda, \gamma_0),R_1(\lambda, \gamma_0))$. From Corollary \ref{pvdev7b-cor1}
we also have  $ F'_{e^{i\theta},\lambda}(z_0)\in V_1\subset\overline{\mathbb{D}}
(C_1(\lambda, \gamma_0),R_1(\lambda, \gamma_0))$ and hence,
$ F'_{e^{i\theta},\lambda}(z_0)\in \partial{V_1}$.

Next, we deal with the uniqueness part. Suppose $ f'(z_0)= F'_{e^{i\theta},\lambda}(z_0)$
for some $f\in\mathcal{S}_1$ and $\theta\in (-\pi, \pi]$. Define
$$h(t)=\frac{\overline{G(z_0)}}{|G(z_0)|}
\left\{f''(z(t))-c_1(z(t),\lambda)\right\}z'(t),
$$
where $\gamma_0: \,z(t)$, $ 0\leq t \leq 1$, as in (\ref{pvdev7b-eq10}).
Then $h(t)$ is a continuous function of $t$ on $[0,1]$ and satisfies the
inequality $|h(t)|\leq r_1(z(t),\lambda)|z'(t)|$. Furthermore, from
(\ref{pvdev7b-eq12}), we get
 \begin{eqnarray*}
 \int_0^1 {\rm Re} \,h(t)\,dt
  & = &
  \int_0^1 {\rm Re} \, \left\{\frac{\overline{G(z_0)}}{|G(z_0)|}
\left\{f''(z(t))-c_1(z(t),\lambda)\right\}z'(t)\right\}dt
\\
 & = &
{\rm Re} \,\left\{\frac{\overline{G(z_0)}}{|G(z_0)|}
\left\{ f'(z_0)-C_1(\lambda, \gamma_0)\right\}\right\}\nonumber
\\
 & = &
 {\rm Re} \,\left\{\frac{\overline{G(z_0)}}{|G(z_0)|}
\left\{ F'_{e^{i\theta},\lambda}(z_0)-C_1(\lambda, \gamma_0)\right\}\right\}\nonumber
\\
 & = &
\int_0^1 r_1(z(t),\lambda)|z'(t)|\,dt\nonumber
 \end{eqnarray*}
which shows that $h(t)= r_1(z(t),\lambda)|z'(t)|$ for all $t\in [0,1]$.
From (\ref{pvdev7b-eq9}) and (\ref{pvdev7b-eq11}) this implies
$f''=F''_{e^{i\theta}, \lambda}$
on $\gamma_0$. From the identity theorem for analytic functions
we have $f''=F''_{e^{i\theta}, \lambda}$
in $\mathbb{D}$ and hence by normalization $f=F_{e^{i\theta},\lambda}$
in $\mathbb{D}$.
\epf

\section{Preparation for the proof of Theorem \ref{pvdev7a-th2}}

\bprop\label{pvdev7a-pro1}
For $f\in\mathcal{S}_2$ with $\lambda\in\mathbb{D}$,
we have
\be\label{pvdev7a-eq4}
\left|f''(z)-c_2(z,\lambda)\right|\leq r_2(z,\lambda), \quad z\in\mathbb{D}
\ee
where
\begin{eqnarray*}
c_2(z,\lambda) & = &
\frac{2M[(1-|z|^2)\lambda+(|z|^2-{|\lambda|}^2)\overline{z}\,]}
{(1-|z|^2)(1+|z|^2-2{\rm Re\,}(\lambda z))},\\
r_2(z,\lambda)& = & \frac{2(1-{|\lambda|}^2)M|z|}
{(1-|z|^2)(1+ |z|^2-2{\rm Re\,}(\lambda z))}.
\end{eqnarray*}
For each $z\in\mathbb{D}\setminus\{0\}$, equality holds if and only if
$f=H_{e^{i\theta},\lambda}$ for some $\theta\in\mathbb{R}$.
\eprop
\bpf
Let $f\in\mathcal{S}_2$. Then (\ref{pvdev7a-eq2}) holds with
$\omega_f(\in\mathcal{{B}}_0)$ and $\omega'_f(0)=\lambda$. It follows from the
Schwarz lemma that
$$\left|\frac{\frac{\omega_f(z)}{z}-\lambda}{1-\overline{\lambda}\frac{\omega_f(z)}{z}}
\right|\leq |z|.
$$
From  (\ref{pvdev7b-eq3a}) and (\ref{pvdev7a-eq2}) we see that this equality is same as
\be\label{pvdev7a-eq6}
\left|\frac{f''(z)-A_2(z,\lambda)}{f''(z)+B_2(z,\lambda)}
\right| \leq |z|\, |\tau_2(z,\lambda)|,
\ee
where
\be\label{pvdev7a-eq7}
A_2(z,\lambda)=\frac{2M\lambda}{1-\lambda z},\quad
B_2(z,\lambda)=\frac{2M}{z-\overline{\lambda}}\quad\mbox { and }
\tau_2(z,\lambda)=\frac{z-\overline{\lambda}}{1-\lambda z}.
\ee
A computation shows that the inequality (\ref{pvdev7a-eq6}) is equivalent to
\be\label{pvdev7a-eq8}
\left|f''(z)-\frac{A_2(z,\lambda)+|z|^2|\tau_2(z,\lambda)|^2 B_2(z,\lambda)}
{1-|z|^2|\tau_2(z,\lambda)|^2}\right| \leq
\frac{|z|\,|\tau_2(z,\lambda)|\,|A_2(z,\lambda)+B_2(z,\lambda)|}{1-|z|^2|\tau_2(z,\lambda)|^2}.
\ee
Also it is easy to obtain that
\begin{eqnarray}
1-|z|^2|\tau_2(z,\lambda)|^2
 & = &
\frac{(1-|z|^2)(1+|z|^2-2{{\rm Re \,} (\lambda z)})}{|1-\lambda z|^2}\nonumber,
\\
A_2(z,\lambda)+ B_2(z,\lambda)
 & = &
\frac{2M(1-{|\lambda|}^2)}{(1-\lambda z)(z-\overline{\lambda})},\nonumber
\end{eqnarray}
and
\begin{eqnarray*}
A_2(z,\lambda)+|z|^2|\tau_2(z,\lambda)|^2 B_2(z,\lambda)
 & = &
\frac{2M[(1-|z|^2)\lambda+(|z|^2-|\lambda|^2)\overline{z}\,]}{|1-\lambda z|^2}\nonumber.
\end{eqnarray*}
Using these,  we obtain that
$$\frac{A_2(z,\lambda)+|z|^2|\tau_2(z,\lambda)|^2 B_2(z,\lambda)}
{1-|z|^2|\tau_2(z,\lambda)|^2}= c_2(z,\lambda)
$$
and
$$\frac{|z|\,|\tau_2(z,\lambda)|\,|A_2(z,\lambda)
+B_2(z,\lambda)|}{1-|z|^2|\tau_2(z,\lambda)|^2}
= r_2(z,\lambda).
$$
Now the inequality (\ref {pvdev7a-eq4}) follows from these equalities
and (\ref{pvdev7a-eq8}). Final part of the proof, namely, the equality case,
follows as in the proof of Proposition \ref{pvdev7b-pro1}.
\epf

The case  $\lambda =0$ of Proposition \ref{pvdev7a-pro1}
gives the following  information.

\bcor\label{pvdev7a-cor3}
Let $f\in \mathcal{S}_2(0)$. Then we have
$$\left|f''(z) -\frac{2M|z|^2\overline{z}}{1-|z|^4}
\right |\leq \frac{2M|z|}{1-|z|^4},
\quad z\in\mathbb{D}.
$$
For each $z\in\mathbb{D}\setminus\{0\}$, equality holds if and only if
$f=H_{e^{i\theta},0}$ for some $\theta\in\mathbb{R}$.
\ecor

In particular, if $f\in \mathcal{S}_2(0)$ then we have
$$(1-|z|^4)\left|f''(z)\right |\leq 2M(1+|z|^2)|z|,
\quad z\in\mathbb{D}
$$
and hence
$$\sup_{z\in\mathbb{D}}(1-|z|^4)\left|f''(z)\right|\leq {4M}.
$$

\bcor\label{pvdev7a-cor1}
Let $\gamma:\,z(t)$, $0\leq t\leq 1$, be a $C^1$-curve in
$\mathbb{D}$ with $z(0)=0$ and $z(1)=z_0$. Then we have
$$
V_2\subset\overline{\mathbb{D}}(C_2(\lambda, \gamma),R_2(\lambda, \gamma))
=\left\{w\in\mathbb{C}:\,|w-C_2(\lambda, \gamma)|\leq R_2(\lambda, \gamma)\right\},
$$
where
$$
C_2(\lambda, \gamma)=\alpha+\int_0^1 c_2(z(t),\lambda)z'(t)\,dt, \quad
R_2(\lambda, \gamma)=\int_0^1 r_2(z(t),\lambda)|z'(t)|\,dt.
$$
\ecor \bpf
The proof is immediate if one uses Proposition \ref{pvdev7a-pro1}
and follows the method of proof of Corollary \ref{pvdev7b-cor1}.
 So we omit the details.
\epf

\blem \cite{samy-vasudev-vuorinen1}\label{pvdev7a-lem2}
For $\theta\in\mathbb{R}$ and $\lambda\in\mathbb{D}$ the function
$$G(z)=\int_0^z \frac{e^{i\theta}\zeta}{\{1+(\overline{\lambda}e^{i\theta}-\lambda)
\zeta-e^{i\theta}{\zeta}^2\}^2} d\zeta ,
 \quad z\in\mathbb{D},
$$
has a double zero at the origin and no zeros elsewhere in $\mathbb{D}$.
Furthermore there exists a starlike univalent function $G_0$ in $\mathbb{D}$
such that $G=2^{-1}e^{i\theta}G^2_0$ and $G_0(0)= G'_0(0)-1=0$.
\elem

\bprop\label{pvdev7a-pro2}
Let $z_0\in\mathbb{D}\setminus \{0\}$. Then for $\theta\in(-\pi,\pi]$
we have $H'_{e^{i\theta},\lambda}(z_0)\in\partial V_2$. Furthermore if
$f'(z_0)= H'_{e^{i\theta},\lambda}(z_0)$ for some
$f\in\mathcal{S}_2$ and $\theta\in(-\pi,\pi]$, then $f= H_{e^{i\theta},\lambda}$.
\eprop
\bpf
From (\ref{pvdev7a_extremal}) we have
$$H''_{a, \lambda}(z)
=\frac{2M\delta (az,\lambda)}{1-\delta(az,\lambda)z}
= \frac{2M(az+\lambda)}{1+(\overline{\lambda}a-\lambda)z-az^2}
$$
and using (\ref{pvdev7a-eq7}), we see that
$$H''_{a, \lambda}(z)-A_2(z,\lambda)
=\frac{2M(1-|\lambda|^2)az}{(1-\lambda z)(1+(\overline{\lambda}a-\lambda)z-az^2)}
$$
and
$$H''_{a, \lambda}(z)+B_2(z,\lambda)
=\frac{2M(1-|\lambda|^2)}{(z-\overline{\lambda}\,)(1+(\overline{\lambda}a-\lambda)z-az^2)}.
$$
Using these, it follows that

\vspace{8pt}
\noindent
$H''_{a, \lambda}(z) -c_2(z,\lambda)$
\begin{eqnarray*}
& = & H''_{a, \lambda}(z)-
\frac{A_2(z,\lambda)+|z|^2|\tau_2(z,\lambda)|^2 B_2(z,\lambda)}
{1-|z|^2|\tau_2(z,\lambda)|^2}\\
&  =  & \frac{1}{1-|z|^2 |\tau_2(z,\lambda)|^2}\left\{H''_{a, \lambda}(z)
-A_2(z,\lambda)-|z|^2|\tau_2(z,\lambda)|^2
\left(H''_{a, \lambda}(z)+B_2(z,\lambda)\right)\right\}\\
& = & \frac{2M(1-|\lambda|^2)az
\{\overline{1+(\overline{\lambda}a-\lambda)z-az^2}\}}
{(1-|z|^2)\{1+|z|^2-2{{\rm Re\,}(\lambda z)}\}\{1+(\overline{\lambda}a-\lambda)z-az^2\}}.
\end{eqnarray*}
Substituting $a=e^{i\theta}$, a computation gives
\begin{eqnarray*}
H''_{e^{i\theta}, \lambda}(z) - c_2(z,\lambda) & = &
\frac{2M(1-|\lambda|^2)e^{i\theta}z
\{\overline{1+(\overline{\lambda}e^{i\theta}-\lambda)z-e^{i\theta}z^2}\}}
{(1-|z|^2)\{1+|z|^2-2{{\rm Re\,}(\lambda z)}\}
\{1+(\overline{\lambda}e^{i\theta}-\lambda)z-e^{i\theta}z^2\}}\\
& = & \frac{2M(1-|\lambda|^2)e^{i\theta}z
\left|1+(\overline{\lambda}e^{i\theta}-\lambda)z-e^{i\theta}z^2\right|^2}
{(1-|z|^2)\{1+|z|^2-2{{\rm Re\,}(\lambda z)}\}
\{1+(\overline{\lambda}e^{i\theta}-\lambda)z-e^{i\theta}z^2\}^2}\\
& = & r_2(z,\lambda)\frac{\left|1+(\overline{\lambda}e^{i\theta}-\lambda)z-e^{i\theta}z^2\right|^2}
{|z|}\frac{e^{i\theta}z}{\{1+(\overline{\lambda}e^{i\theta}-\lambda)z-e^{i\theta}z^2\}^2}.
\end{eqnarray*}
Using Lemma \ref{pvdev7a-lem2}, we may rewrite the last equality as
\be\label{pvdev7a-eq9}
H''_{e^{i\theta}, \lambda}(z)-c_2(z,\lambda)
= r_2(z,\lambda)\frac{G'(z)}{|G'(z)|}
\ee
where $G(z)$ is defined as in Lemma \ref{pvdev7a-lem2}.
Since $G_0$ defined by  $G=2^{-1}e^{i\theta}G^2_0$
is a normalized starlike function, for any $z_0\in\mathbb{D}\setminus\{0\}$, the line
segment joining $0$ and $G_0(z_0)$ entirely lies in $G_0(\mathbb{D})$.
As before, define $\gamma_0$ by
$$\gamma_0:\,z(t)=G_0^{-1}(tG_0(z_0)), \quad 0\leq t \leq 1.
$$
We observe that $G(z(t))= 2^{-1}e^{i\theta}G_0(z(t))^2= 2^{-1}e^{i\theta}(tG_0(z_0))^2
=t^2 G(z_0)$ and so, we get
\be\label{pvdev7a-eq11}
G'(z(t))z'(t)=2tG(z_0), \quad t\in [0,1].
\ee
From this, (\ref{pvdev7a-eq9}) and proceeding exactly as in the proof of
Proposition \ref{pvdev7b-pro2}, we end up with
\be\label{pvdev7a-eq12}
H'_{e^{i\theta},\lambda}(z_0)-C_2(\lambda, \gamma_0)
=\frac{G(z_0)}{|G(z_0)|}R_2(\lambda, \gamma_0)
\ee
which gives $H'_{e^{i\theta},\lambda}(z_0)\in\partial{\overline{\mathbb{D}}}
(C_2(\lambda, \gamma_0),R_2(\lambda, \gamma_0))$.
From Corollary \ref{pvdev7a-cor1}, we also have
$ H'_{e^{i\theta},\lambda}(z_0)\in V_2\subset\overline{\mathbb{D}}
(C_2(\lambda, \gamma_0),R_2(\lambda, \gamma_0))$.
Hence, $ H'_{e^{i\theta},\lambda}(z_0)\in \partial V_2$.

Uniqueness part follows similarly. Indeed, suppose $ f'(z_0)= H'_{e^{i\theta},\lambda}(z_0)$
for some $f\in\mathcal{S}_2$ and $\theta\in (-\pi, \pi]$ and introduce,
$$h(t)=\frac{\overline{G(z_0)}}{|G(z_0)|}\left\{f''(z(t))-c_2(z(t),\lambda)\right\}z'(t).
$$
Then $h(t)$ is continuous function of $t\in [0,1]$ and satisfies
$|h(t)|\leq r_2(z(t),\lambda)|z'(t)|$. Furthermore, from (\ref{pvdev7a-eq12}), we
easily see that
\begin{eqnarray*}
\int_0^1 {\rm Re} \,h(t)\,dt& = &
\int_0^1 r_2(z(t),\lambda)|z'(t)|\,dt .
\end{eqnarray*}
Thus, $h(t)= r_2(z(t),\lambda)|z'(t)|$ for all $t\in [0,1]$.
From (\ref{pvdev7a-eq9}) and (\ref{pvdev7a-eq11}) this implies that
$f''=H''_{e^{i\theta}, \lambda}$ on $\gamma_0$. From the identity theorem for analytic
functions we deduce that $f''=H''_{e^{i\theta}, \lambda}$ in $\mathbb{D}$ and hence
by normalization $f=H_{e^{i\theta},\lambda}$ in $\mathbb{D}$.
\epf
\section{Proofs of Theorems \ref{pvdev7b-th1} and \ref{pvdev7a-th2}}
\begin{nonsec}
Proof of Theorem \ref{pvdev7b-th1}
\end{nonsec}  We prove that the closed curve
$(-\pi,\pi]\ni \theta \mapsto F'_{e^{i\theta},\lambda}(z_0)$
is simple. Suppose that $F'_{e^{i\theta_1},\lambda}(z_0)=
F'_{e^{i\theta_2},\lambda}(z_0)$ for some
$\theta_1,\theta_2\in(-\pi,\pi]$ with $\theta_1\neq\theta_2$.
Then, from Proposition \ref{pvdev7b-pro2}, we conclude that
$F_{e^{i\theta_1},\lambda}= F_{e^{i\theta_2},\lambda}$.

From (\ref{pvdev7b-eq3}) and (\ref{pvdev7b-eq7}) we have
$$
\tau_1\left(\frac{\omega_{F_{e^{i\theta},\lambda}}}{z},\lambda\right)
=\frac{(\overline{\beta}+\overline{\lambda}^2)e^{i\theta}z+(\overline{\beta}\lambda+\overline{\lambda})}
{(\overline{\lambda}+\overline{\beta}\lambda)e^{i\theta}z+(\overline{\beta}{\lambda}^2+1)}.
$$
Since $F_{e^{i\theta_1},\lambda}= F_{e^{i\theta_2},\lambda}$, we have the following relation
$$
\tau_1\left(\frac{\omega_{F_{e^{i\theta_1},\lambda}}}{z},\lambda\right)
=\tau_1\left(\frac{\omega_{F_{e^{i\theta_2},\lambda}}}{z},\lambda\right).
$$
That is
$$
\frac{(\overline{\beta}+\overline{\lambda}^2)e^{i\theta_1}z+(\overline{\beta}\lambda+\overline{\lambda})}
{(\overline{\lambda}+\overline{\beta}\lambda)e^{i\theta_1}z+(\overline{\beta}{\lambda}^2+1)}
=\frac{(\overline{\beta}+\overline{\lambda}^2)e^{i\theta_2}z+(\overline{\beta}\lambda+\overline{\lambda})}
{(\overline{\lambda}+\overline{\beta}\lambda)e^{i\theta_2}z+(\overline{\beta}{\lambda}^2+1)}.
$$
By a simplification, the last expression implies
$$e^{i\theta_1}z=e^{i\theta_2}z
$$
which is a contradiction for the choice of
$\theta_1$ and $\theta_2$. Thus the curve is simple.

Since $V_1$ is a compact convex subset of $\mathbb{C}$ and has nonempty interior,
the boundary $\partial{V_1}$ contains the curve
$(-\pi,\pi]\ni \theta \mapsto F'_{e^{i\theta},\lambda}(z_0)$. Note that a simple
closed curve cannot contain any simple closed curve other than itself. Thus
$\partial{V_1}$ is given by
$(-\pi,\pi]\ni \theta \mapsto  F'_{e^{i\theta},\lambda}(z_0)$.
\begin{nonsec}
Proof of Theorem \ref{pvdev7a-th2}
\end{nonsec}
We prove that the closed curve
$(-\pi,\pi]\ni \theta \mapsto H'_{e^{i\theta},\lambda}(z_0)$
is simple. Suppose that $H'_{e^{i\theta_1},\lambda}(z_0)=
H'_{e^{i\theta_2},\lambda}(z_0)$ for some
$\theta_1,\theta_2\in(-\pi,\pi]$ with $\theta_1\neq\theta_2$. Then
from Proposition \ref{pvdev7a-pro2} we have $
H_{e^{i\theta_1},\lambda}= H_{e^{i\theta_2},\lambda}$. From
(\ref{pvdev7a-eq3}) and (\ref{pvdev7a-eq7}) this shows a contradiction
$$
e^{i\theta_1}z=\tau_2\left(\frac{\omega_{H_{e^{i\theta_1},\lambda}}}{z},\lambda\right)
=\tau_2\left(\frac{\omega_{H_{e^{i\theta_2},\lambda}}}{z},\lambda\right)=e^{i\theta_2}z.
$$
Thus the curve is simple.

Again, since $V_2$ is a compact convex subset of $\mathbb{C}$ and has nonempty interior,
the boundary $\partial V_2$ contains the curve
$(-\pi,\pi]\ni \theta \mapsto H'_{e^{i\theta},\lambda}(z_0)$.
The same reasoning as in the proof of Theorem \ref{pvdev7b-th1} shows that
$\partial V_2$ is given by
$(-\pi,\pi]\ni \theta \mapsto  H'_{e^{i\theta},\lambda}(z_0)$.

\section{Geometric view of Theorems \ref{pvdev7b-th1} and \ref{pvdev7a-th2}}

  Using Mathematica 4.1,  we describe the boundary of the sets
$V_j(z_0, \lambda)$ for $j=1,2$.  Here we give the Mathematica program
which is used to plot the boundary of the sets  $V_j(z_0, \lambda)$ for $j=1,2$.
We refer \cite{Ruskeepaa} for Mathematica program. The short notations in this program
are of the form: ``z0 for $z_0$'', ``a for $\alpha$'', ``lam for $\lambda$'',
``m for $M$'' and ``b for $\beta$''.

\vspace{1cm}

{\tt
\begin{verbatim}
Remove["Global`*"];
(* The values ``z0, a, lam, m, b'' are for FIGURE 1 *)
z0 = 0.00882581 - 0.514124I
a = -230.939 + 799.526I
lam = 0.427174 + 0.0755107I
m = 509.317
b = 0.94485 + 0.0416585I

Q1[b_, m_, lam_, the_] :=
m((Exp[I*the]z + lam)z +b(1 + Conjugate[lam]Exp[I*the]z))/
((1 +(Conjugate[lam]*Exp[I*the] +Conjugate[b]*lam)*z)+
Conjugate[b]*Exp[I*the]*z*z);


myf1[a_, b_, m_, lam_, the_, z0_] :=
a +NIntegrate[Q1[b, m, lam, the], {z, 0, z0}];

image1 = ParametricPlot[{Re[myf1[a, b, m, lam, the, z0]],
        Im[myf1[a, b, m, lam, the, z0]]}, {the, -Pi, Pi},
        AspectRatio -> Automatic];


Clear[a, b, m, lam, the, z0, myf1];

z0 = 0.00882581 - 0.514124I
a = -230.939 + 799.526I
lam = 0.839567
m = 0.254877

Q2[m_, lam_, the_] :=
2*m*(Exp[I*the]*z + lam)/(1 + lam*(Exp[I*the] - 1)*z
- Exp[I*the]*z*z);

myf2[a_, m_, lam_, the_, z0_] :=
a + NIntegrate[Q2[m, lam, the], {z, 0, z0}];

image2 = ParametricPlot[{Re[myf2[a, m, lam, the, z0]],
Im[myf2[a, m, lam, the, z0]]}, {the, -Pi, Pi},
AspectRatio -> Automatic];

image=Show[GraphicsArray[{image1,image2},GraphicsSpacing --->0.5]]

Clear[a, m, lam, the, z0, myf2];

\end{verbatim}
}

\vspace{0.5cm}

The following pictures give the geometric view of the boundary of the
sets $V_j(z_0, \lambda)$ for each $j=1,2$.
In each of the following figures the left hand side figure describes the boundary of the
set $V_1(z_0, \lambda)$ for each fixed value of $z_0\in\mathbb{D}\setminus\{0\}$, $\lambda\in\mathbb{D}$,
$\alpha, \beta\in\mathbb{C}$ and $M\in\mathbb{R}$  with  $0<M\leq |\alpha|$
and $|\beta|\leq 1$. These values are given in the  first column of each the figure.
Similarly the right hand side of each of the following figures describes the boundary
of the set $V_2(z_0, \lambda)$ for each fixed value of $z_0\in\mathbb{D}\setminus\{0\}$,
$\lambda\in[0,1)$, $\alpha\in\mathbb{C}\setminus \{0\}$ and $M\in\mathbb{R}$ such that $0<M\leq 1/{\log 4}$
and these values are given in the second column of each of the figures.
Note that according to  Proposition \ref{pvdev7b-pro3} and Proposition \ref{pvdev7a-pro3}
the regions bounded by the curves $\partial V_j(z_0, \lambda)$ for  $j=1,2$ are compact and convex.

\begin{figure}[htp]
\begin{center}
\includegraphics[width=5cm]{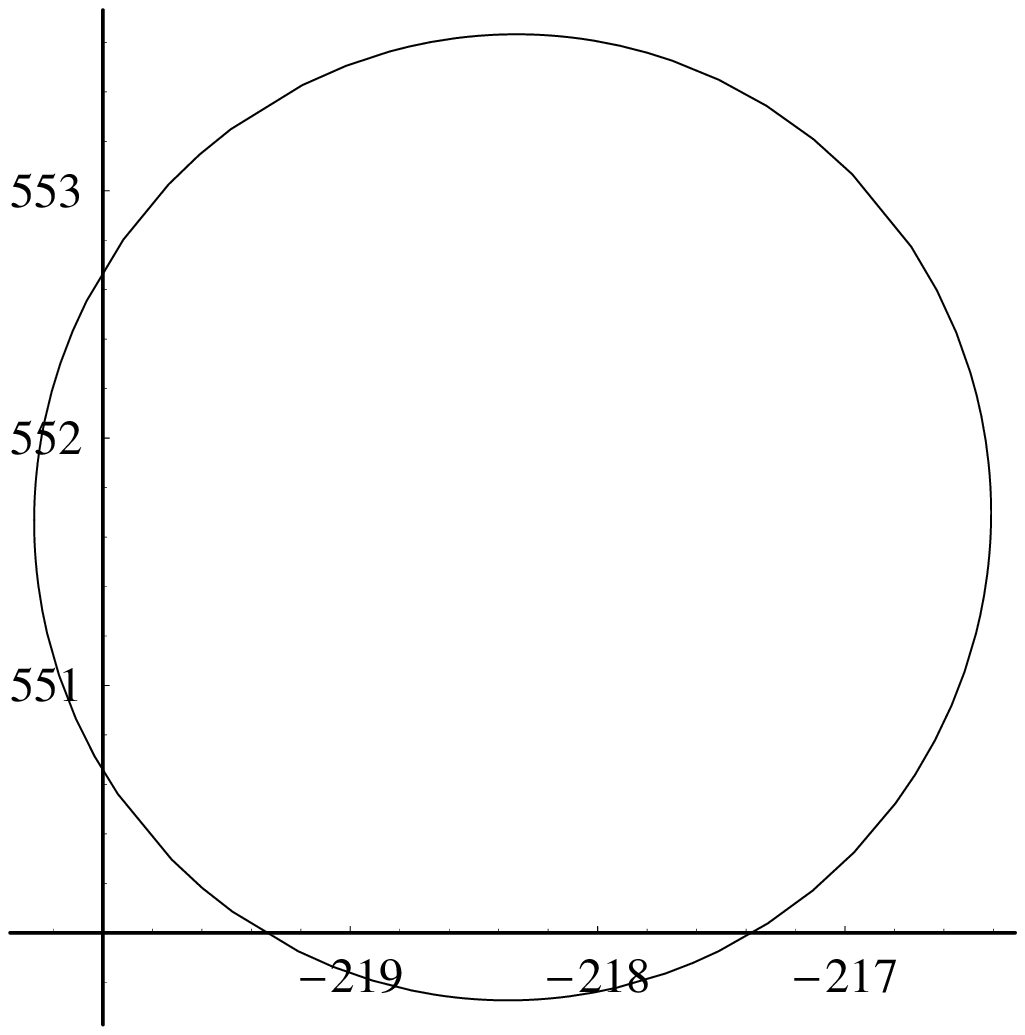}
\hspace{1cm}
\includegraphics[width=5cm]{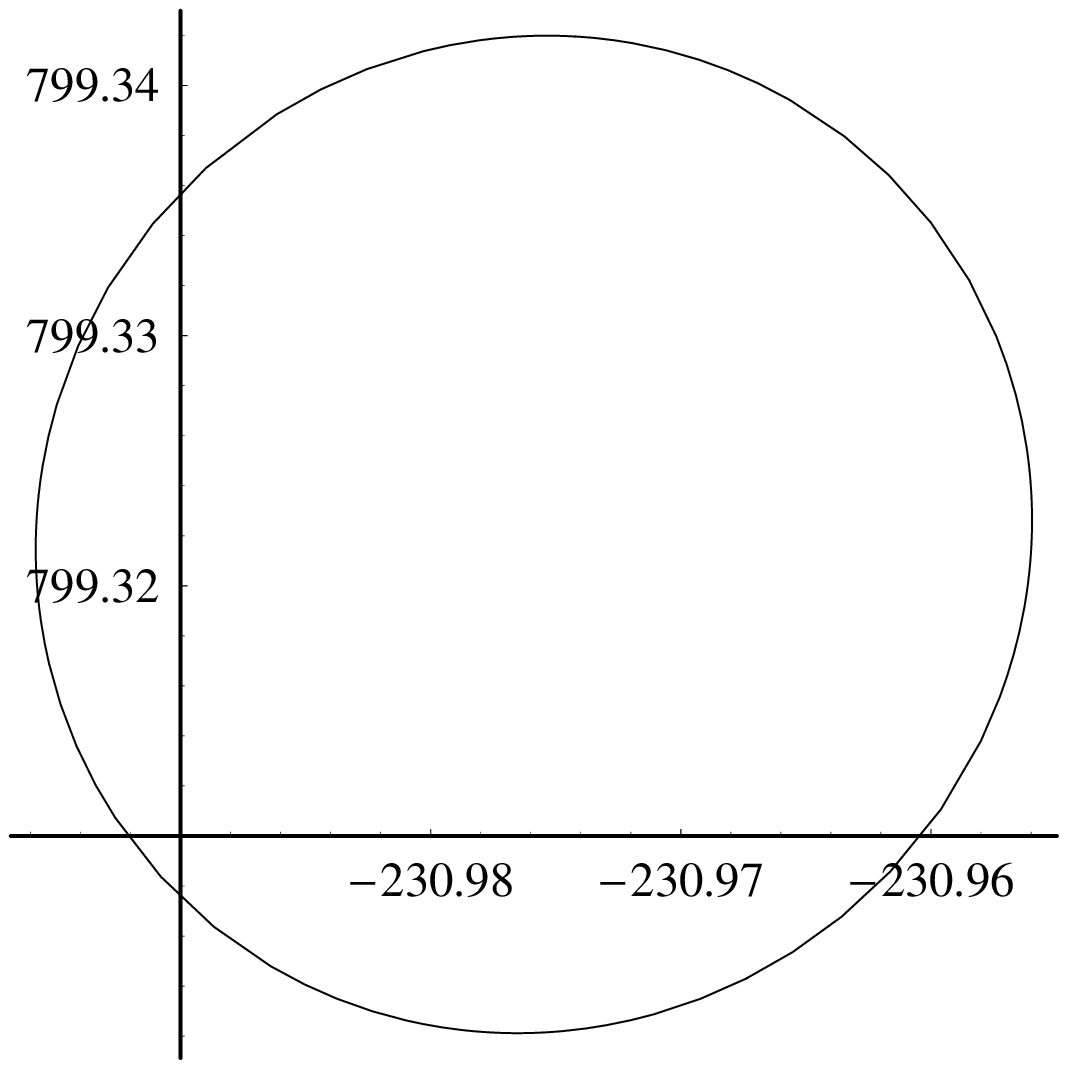}
\end{center}
$\partial V_1(z_0,\lambda)$ \hspace{5cm}  $\partial V_2(z_0,\lambda)$
\caption{Region of variability for $f'(z_0)$}
\end{figure}
\begin{center}
$\begin{array}{ll}
z_0 =  0.00882581 - 0.514124i       \hspace{3cm} & z_0 =   0.00882581 - 0.514124i \\
\alpha =-230.939 + 799.526i                & \alpha = -230.939 + 799.526i \\
\lambda =0.427174 + 0.0755107i                    & \lambda =  0.839567\\
M = 509.317                     & M = 0.254877\\
\beta = 0.94485 + 0.0416585i         &
\end{array}$
\end{center}

\begin{figure}[htp]
\begin{center}
\includegraphics[width=5cm]{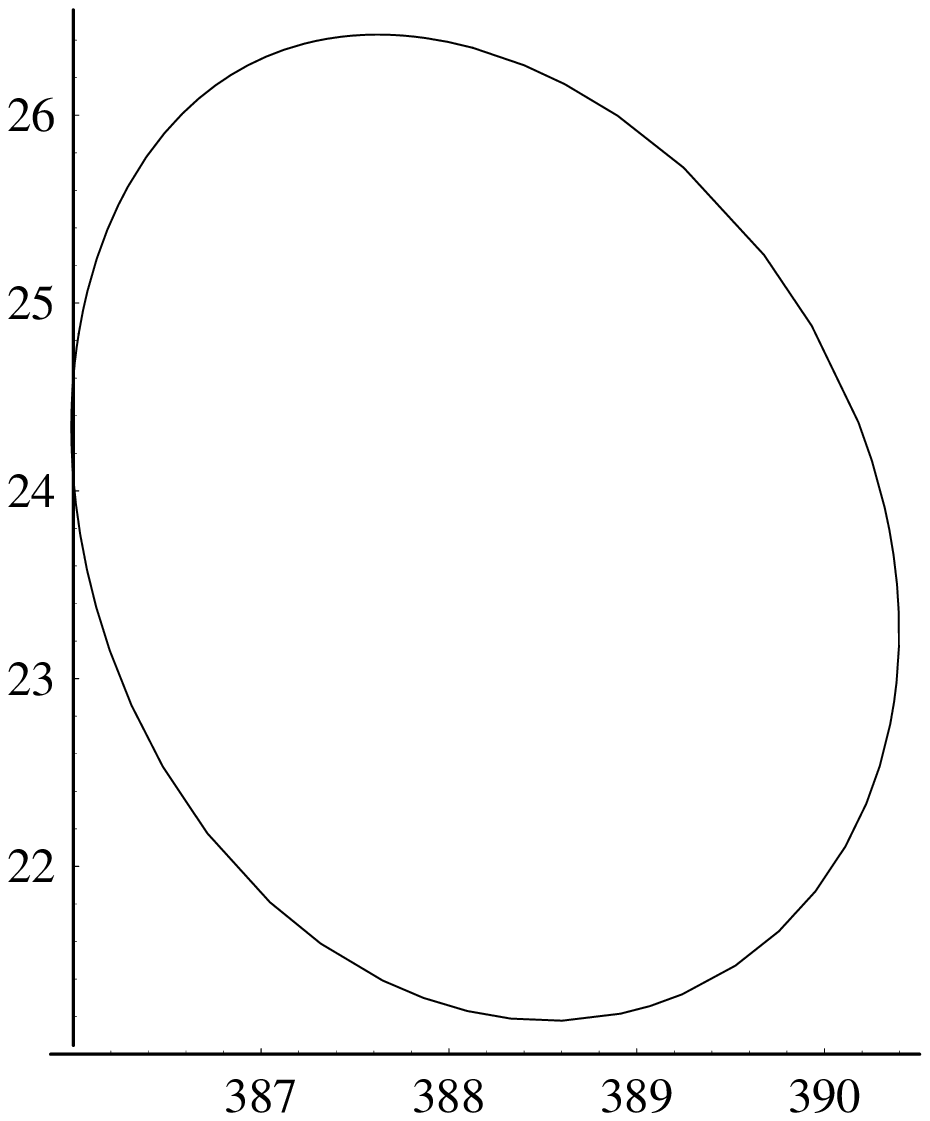}
\hspace{1cm}
\includegraphics[width=6cm]{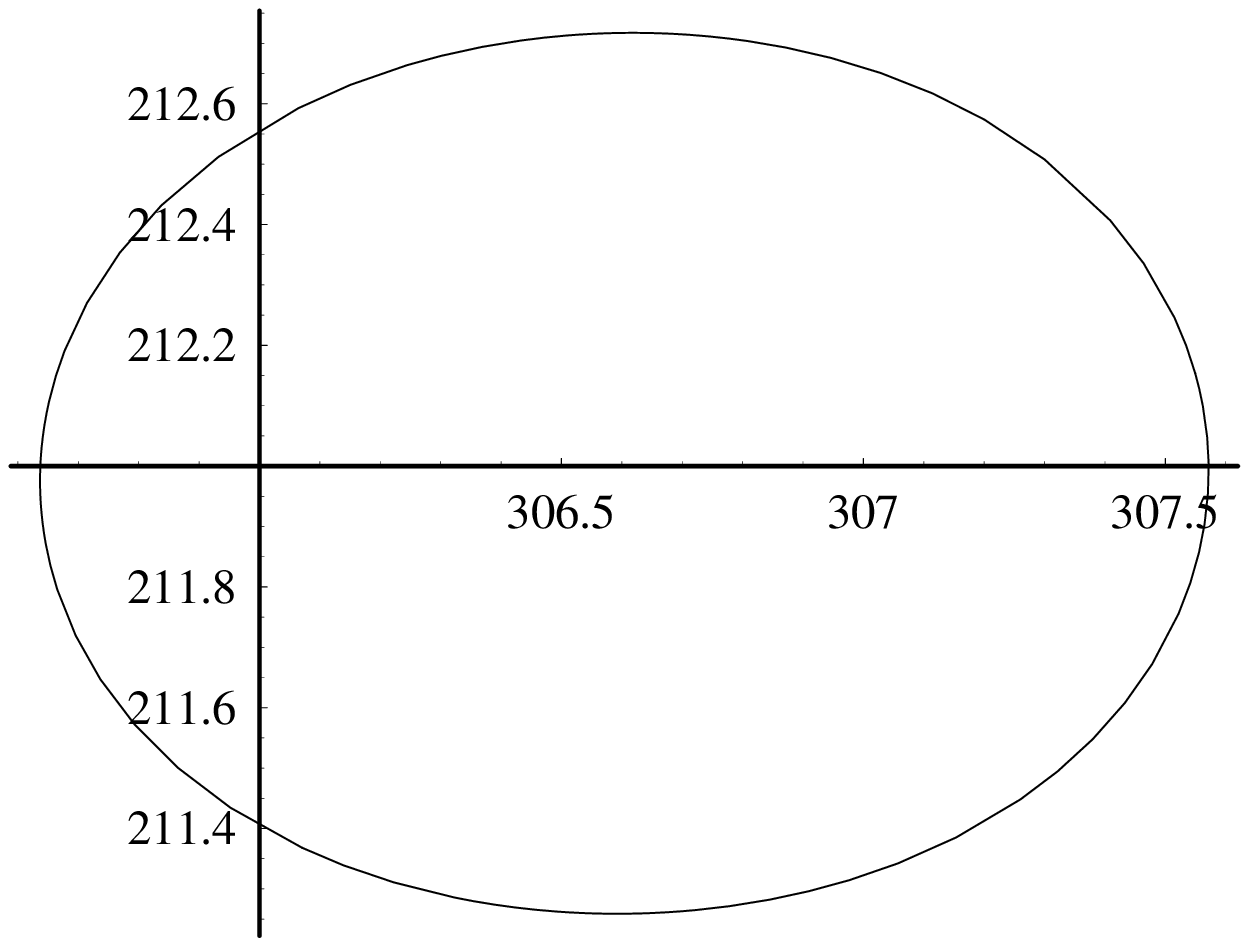}
\end{center}
$\partial V_1(z_0,\lambda)$ \hspace{5cm}  $\partial V_2(z_0,\lambda)$
\caption{Region of variability for $f'(z_0)$}
\end{figure}
\begin{center}
$\begin{array}{ll}
z_0 =-0.439619 - 0.843107i      \hspace{3cm} & z_0 = -0.439619 - 0.843107i\\
\alpha =306.095 + 212.047i             & \alpha = 306.095 + 212.047i \\
\lambda =-0.847689 - 0.07592i                 & \lambda =    0.0802624\\
M = 206.329                    & M = 0.673609\\
\beta = 0.67079 + 0.843107i                    &
\end{array}$
\end{center}

\begin{figure}[htp]
\begin{center}
\includegraphics[width=5cm]{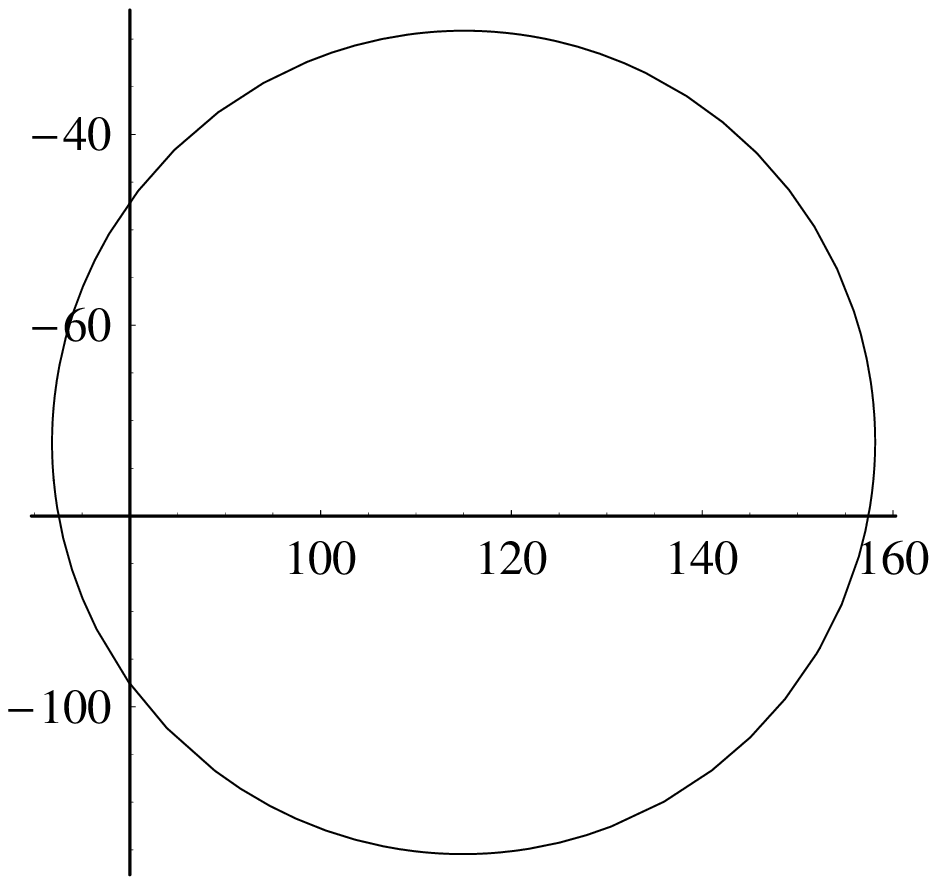}
\hspace{1cm}
\includegraphics[width=7cm]{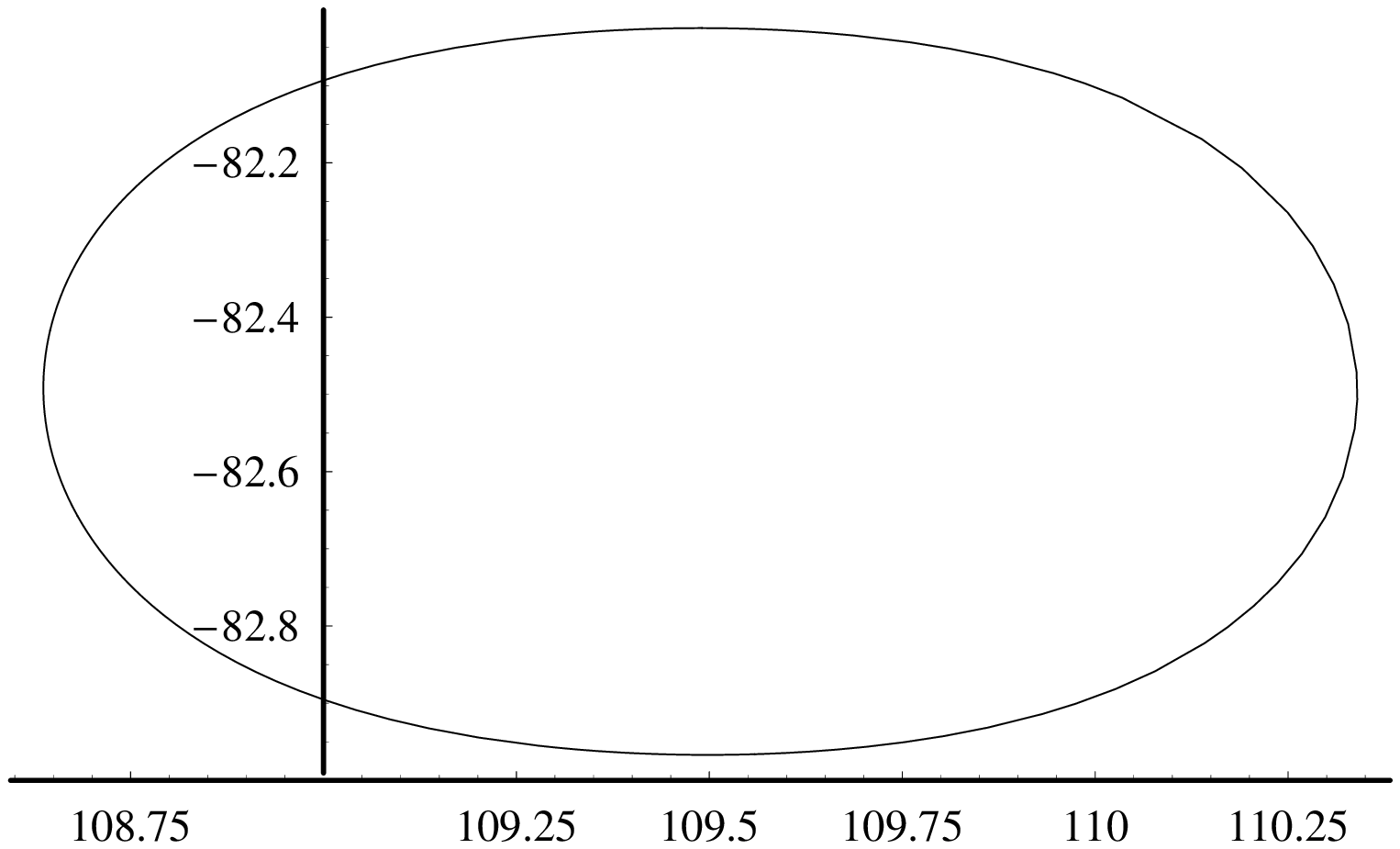}
\end{center}
$\partial V_1(z_0,\lambda)$ \hspace{5cm}  $\partial V_2(z_0,\lambda)$
\caption{Region of variability for $f'(z_0)$}
\end{figure}
\begin{center}
$\begin{array}{ll}
z_0 =-0.971007 + 0.211382i       \hspace{3cm} & z_0 = -0.971007 + 0.211382i\\
\alpha =108.958 - 82.5096i              & \alpha = 108.958 - 82.5096i \\
\lambda =0.0327389 - 0.0219389i                 & \lambda =   0.148939\\
M = 132.988                     & M = 0.390188\\
\beta = -0.0264629 - 0.114565i                    &
\end{array}$
\end{center}

\begin{figure}[htp]
\begin{center}
\includegraphics[width=5cm]{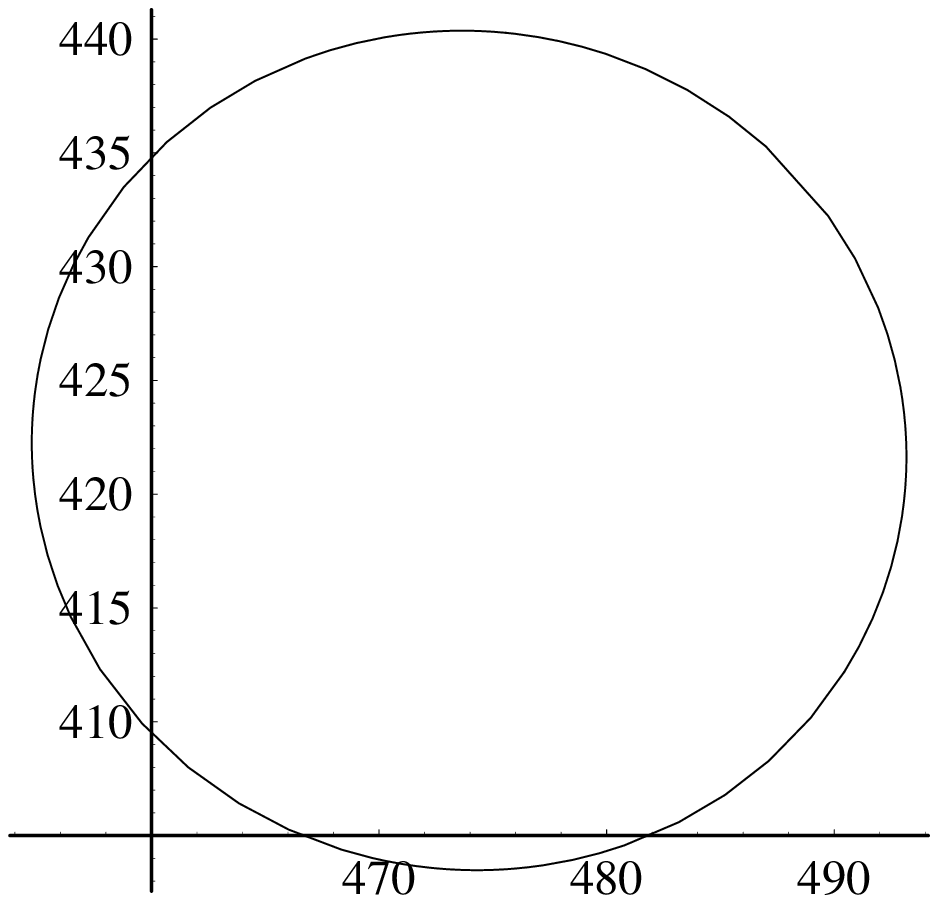}
\hspace{1cm}
\includegraphics[width=7cm]{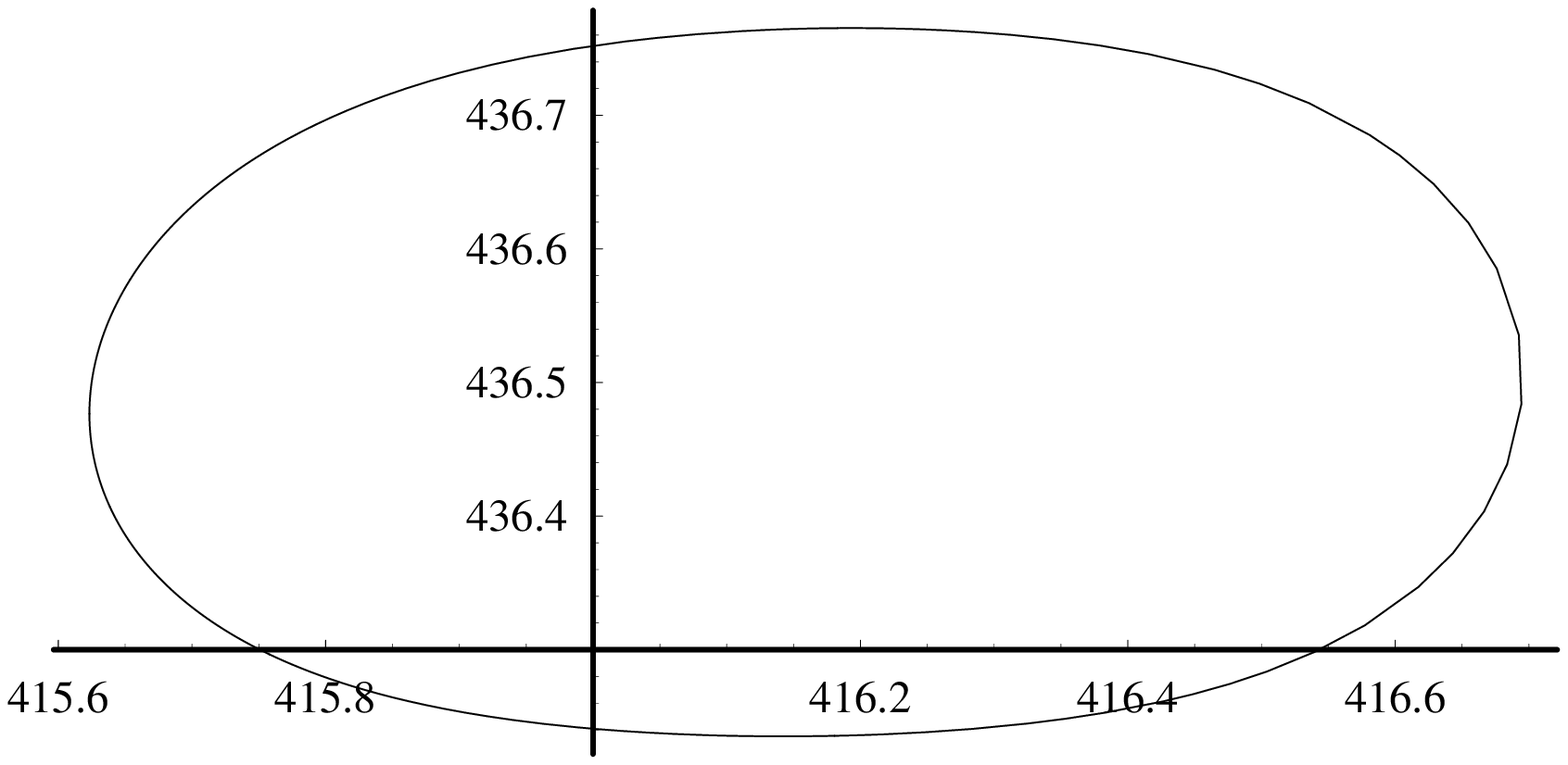}
\end{center}
$\partial V_1(z_0,\lambda)$ \hspace{5cm}  $\partial V_2(z_0,\lambda)$
\caption{Region of variability for $f'(z_0)$}
\end{figure}
\begin{center}
$\begin{array}{ll}
z_0 =-0.844358 - 0.529996i       \hspace{3cm} & z_0 = -0.844358 - 0.529996i \\
\alpha =416.349 + 436.752i               & \alpha = 416.349 + 436.752i \\
\lambda =-0.0872118 + 0.664418i                  & \lambda =  0.7262\\
M = 97.2626                      & M = 0.620559\\
\beta = -0.549327 + 0.592394i                    &
\end{array}$
\end{center}

\newpage

\begin{figure}[htp]
\begin{center}
\includegraphics[width=5cm]{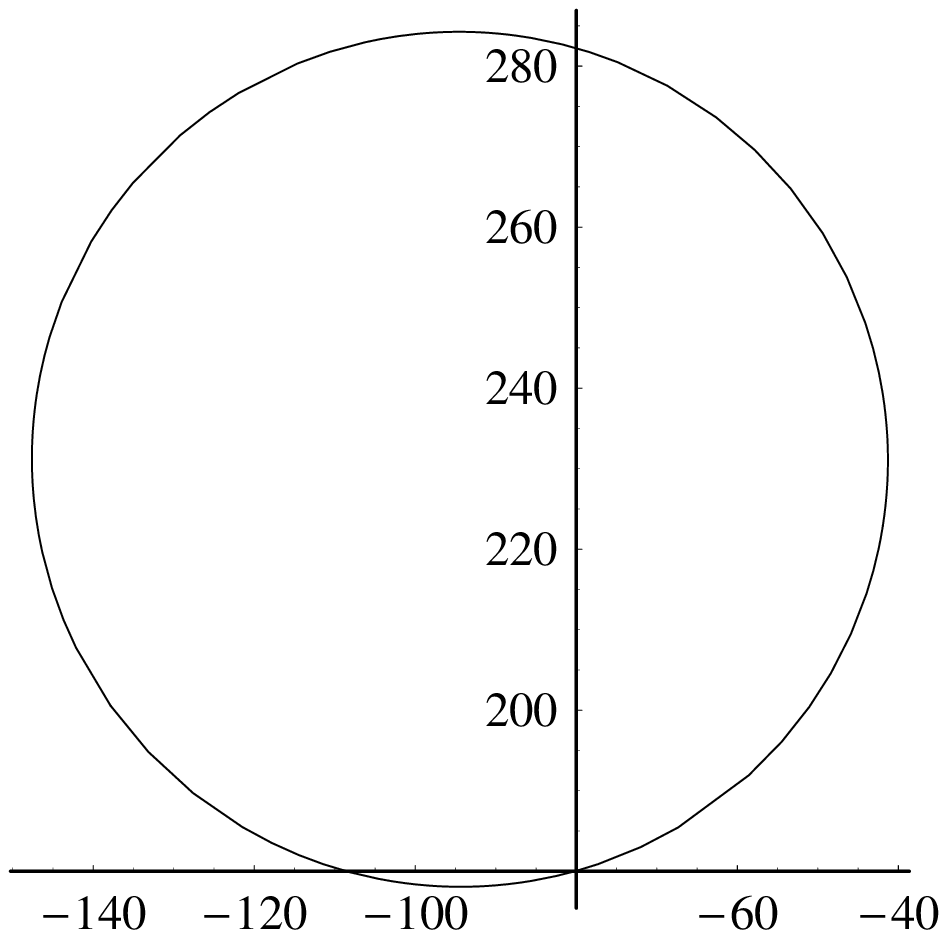}
\hspace{1cm}
\includegraphics[width=7cm]{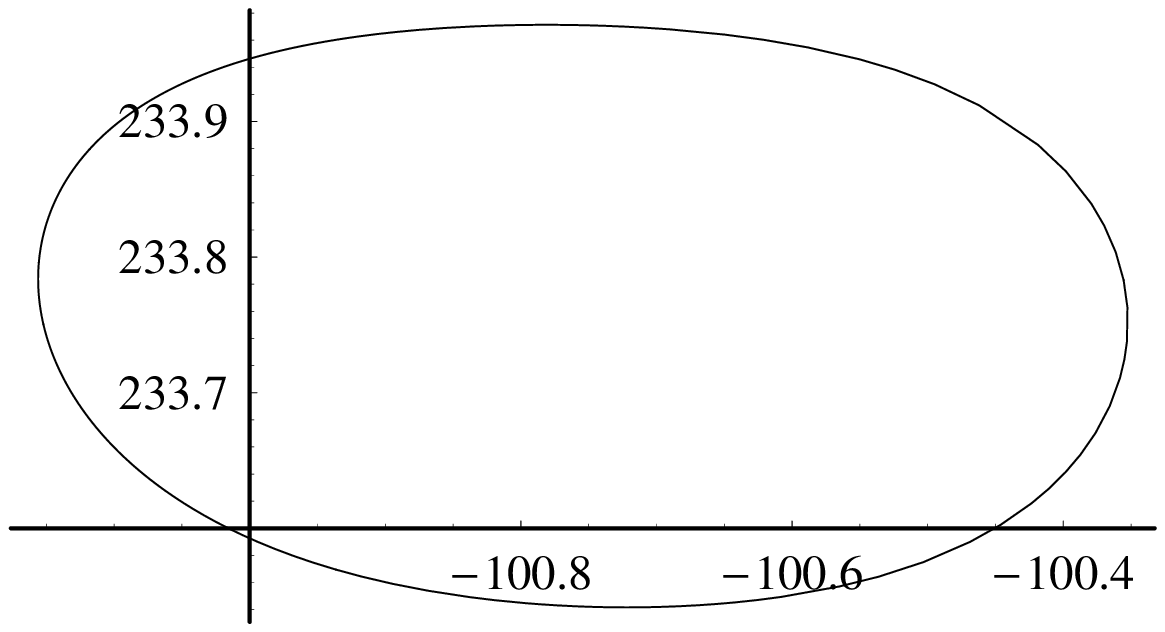}
\end{center}
$\partial V_1(z_0,\lambda)$ \hspace{5cm}  $\partial V_2(z_0,\lambda)$
\caption{Region of variability for $f'(z_0)$}
\end{figure}
\begin{center}
$\begin{array}{ll}
z_0 =-0.605185 + 0.789592i       \hspace{3cm} & z_0 = -0.605185 + 0.789592i \\
\alpha =-100.796 + 233.556i                & \alpha = -100.796 + 233.556i \\
\lambda =0.0523661 + 0.167249i                   & \lambda =  0.63945\\
M = 164.079                      & M = 0.354197\\
\beta = 0.00810121 - 0.00819085i                    &
\end{array}$
\end{center}

\end{document}